\documentclass[11pt]{article}
\date{}
\makeatletter
\@addtoreset{equation}{section}
\makeatother
\usepackage{amssymb,amsmath}
\usepackage{cite}
\usepackage{hyperref}

\begin{document}

\title{\bf Normalized Solutions to Nonautonomous Kirchhoff Equation \footnote{Supported by National Natural Science Foundation of China(No. 12271443) and Natural Science Foundation
of Chongqing, China(cstc2020jcyjjqX0029).}}
\author{{Xin Qiu,\ Zeng-Qi Ou,\ Ying Lv \footnote{Corresponding author.
E-mail address: qxinqiuxin@163.com(X. Qiu), ouzengq707@sina.com(Z.-Q. Ou), ly0904@swu.edu.cn(Y. Lv).}}\\
{\small \emph{ School  of  Mathematics  and  Statistics, Southwest University,  Chongqing {\rm400715},}}\\
{\small \emph{People's Republic of China}}}
\maketitle
\baselineskip 17pt

{\bf Abstract:} In this paper, we study the existence of normalized solutions to the following Kirchhoff equation with a perturbation:
$$
\left\{
\begin{aligned}
&-\left(a+b\int _{\mathbb{R}^{N}}\left | \nabla u \right|^{2} dx\right)\Delta u+\lambda u=|u|^{p-2} u+h(x)\left |u\right |^{q-2}u, \quad  \text{ in } \mathbb{R}^{N}, \\
&\int_{\mathbb{R}^{N}}\left|u\right|^{2}dx=c, \quad u \in H^{1}(\mathbb{R}^{N}),
\end{aligned}
\right.
$$
where $1\le N\le 3, a,b,c>0, 1\leq q<2$, $\lambda \in \mathbb{R}$. We treat three cases.\\
(i)When $2<p<2+\frac{4}{N}, h(x)\ge0$, we obtain the existence of global constraint minimizers.\\
(ii)When $2+\frac{8}{N}<p<2^{*}, h(x)\ge0$, we prove the existence of mountain pass solution.\\
(iii)When $2+\frac{8}{N}<p<2^{*}, h(x)\leq0$, we establish the existence of bound state solutions.

  \par

\section{Introduction and main results}

\ \ \ \ \ \ In this paper, we consider the existence of solutions with prescribed $L^{2}$-norm to the following Kirchhoff problem with a perturbation
\begin{equation}\label{initial}
\left\{
\begin{aligned}
&-\left( a+b\int_{\mathbb{R}^{N}}\left| \nabla u \right |^{2}dx \right )\Delta u+\lambda u=|u|^{p-2}u+h(x)\left|u\right|^{q-2}u,\quad \text { in } \mathbb{R}^{N}, \\
&\int_{\mathbb{R}^{N}}\left|u\right|^{2}dx=c, \quad u \in H^{1}(\mathbb{R}^{N}),
\end{aligned}
\right.
\end{equation}
where $1\le N\le 3, a,b,c>0, p\in (2,2^{*}), q \in [1,2)$, $h(x):\mathbb{R}^{N}\to \mathbb{R}$ is a potential, $2^{*}=6$ if $N=3$, and $2^{*}=+\infty$ if $N=1,2$. Based on these observations, we establish the existence of normalized solutions under different assumptions on $h(x)$.

The energy functional of Eq.\eqref{initial} is defined by
\begin{equation}\label{Iu}
I(u)=\frac{a}{2} \int_{\mathbb{R}^{N}}\left|\nabla u\right|^{2}dx+\frac{b}{4} \left (\int_{\mathbb{R}^{N}}\left|\nabla u\right|^{2}dx \right )^{2}-\frac{1}{p} \int_{\mathbb{R}^{N}}\left|u\right|^{p} dx-\frac{1}{q} \int_{\mathbb{R}^{N}} h(x)\left|u\right|^{q}dx
\end{equation}
constrained on the $L^{2}$-spheres in $H^{1}(\mathbb{R}^{N})$:
$$
S_{c}=\left\{u\in H^{1}(\mathbb{R}^{N}): \|u\|_{2}^{2}=c>0\right\}.
$$

In 1883, Kirchhoff \cite{z1} first proposed the following nonlinear wave equation 
$$
\rho\frac{\partial^{2}u}{\partial t^{2}}-\left(\frac{P_{0}}{h}+\frac{E}{2L}\int_{0}^{L}\left|\frac{\partial u}{\partial x}\right|^{2}dx\right)\frac{\partial^{2}u}{\partial x^{2}}=0,
$$
which extends the original wave equation by describing the transversal oscillations of a stretched string and particularly considering the subsequent change in string length caused by oscillations. Thereafter, there has been a boom in the study of the Kirchhoff-type equation. We can refer to \cite{A1,C2,D1} for the physical background about Kirchhoff problem.


Mathematically, Eq.\eqref{initial} is not a pointwise identity as a result of the emergence of the term $\left(b\int_{\mathbb{R}^{N}}\left|\nabla u\right|^{2}dx\right)\Delta u$. This causes some mathematical difficulties. In the renowned paper \cite{L1}, J.L. Lions raised an abstract framework which has received much attention. There are two ways to study the Kirchhoff-type equation. The first approach is to consider fixing the parameter $\lambda \in \mathbb{R}$. In this case, there are a lot of results which have been widely studied by using variational methods. We can refer to \cite{F1,G1,M1,M2,H1} and the references therein. Another way is to fix the $L^{2}$-norm. In this case, the desired solutions have a priori prescribed $L^{2}$-norm, which are usually referred to as normalized solutions in the literature, that is, for any fixed $c>0$, we take $\left(u_{c}, \lambda_{c}\right)\in H^{1}(\mathbb{R}^{N})\times \mathbb{R}$ as normalized solution with $\|u_{c}\|_{2}^{2}=c$, $\lambda_{c}$ is a Lagrange multiplier. From a physical perspective, this type of the prescribed mass problem has physically significant in Bose-Einstein condensates and the nonlinear optics framework.

For the local case, i.e., $b=0$, Eq.\eqref{initial} reduces to the general Schr\"{o}dinger type:
\begin{equation}\label{yiban}
\left\{\begin{aligned}
&-\Delta u+\lambda u=f(x,u), \quad \text { in } \mathbb{R}^{N}, \\
&\int_{\mathbb{R}^{N}} \left|u\right|^{2} dx=c, \quad u \in H^{1}(\mathbb{R}^{N}),
\end{aligned}\right.
\end{equation}
which dates back to the groundbreaking work by Stuart. In \cite{S1,S2}, Stuart tackled problem \eqref{yiban} for $f(x,u)=|u|^{p-2}u$ and $p\in(2,2+\frac{4}{N})$ ($L^2$-subcritical case), here $2+\frac{4}{N}$ is called the $L^2$-critical exponent. For $L^2$-subcritical case, minimization method is the conventional method to find normalized solutions. When $f$ is $L^2$-supercritical growth, one of the groundbreaking work in the $L^{2}$-supercritical case is accomplished by Jeanjean\cite{r4}. Jeanjean developed a novel argument related to the mountain pass geometry by the the stretched functional. Bartsch and Soave \cite{soave1,soave2} also proposed a new approach by using a minimax principle based on the homotopy stable family to prove the existence of normalized solutions for problem \eqref{yiban}. Moreover, Soave in \cite{Z2} has been studied the combined nonlinearity case $f(x,u)=|u|^{p-2}u+\mu |u|^{q-2}u$, $2<q\leq 2+\frac{4}{N}\leq p<2^*$ and $q<p$, where $2^*=\infty$ if $N\leq2$ and $2^*=\frac{2N}{N-2}$ if $N\geq3$. Soave showed that nonlinear terms with different power strongly affects the geometry of the functional and the existence and properties of ground states.

When $f(x,u)=a(x)f(u)$, the solutions to the nonautonomous problem
which studied by Chen and Tang \cite{C4} firstly. Compared with the autonomous
problems, the main challenge of the problem 
is constructing a $(PS)$ sequence with additional property to recover the compactness. Very recently, Chen and Zou \cite{z18} studied the following problem with a perturbation
\begin{equation}\label{ckeq}
\left\{\begin{aligned}
&-\Delta u+\lambda u=|u|^{p-2}u+h(x), \quad \text { in } \mathbb{R}^{N}, \\
&\int_{\mathbb{R}^{N}} \left|u\right|^{2} dx=c, \quad u \in H^{1}(\mathbb{R}^{N}),
\end{aligned}\right.
\end{equation}
where $h(x)\ge 0$. For $p\in(2,2+\frac{4}{N})$ and arbitrarily positive perturbation, Chen and Zou proved that there exists a global minimizer with negative energy. The existence of a mountain pass solution with positive energy for $p\in(2+\frac{4}{N},2^{*})$ has been studied. We can see \cite{A2,C5,Z1} for more details.

For the nonlocal case, i.e., $b>0$, the more general form of Eq.\eqref{initial} is the following equation
\begin{equation}\label{yiban4}
\left\{\begin{aligned}
&-\left(a+b\int_{\mathbb{R}^{N}}\left | \nabla u \right |^{2}dx\right)\Delta u+\lambda u=f(x,u), \quad \text { in } \mathbb{R}^{N}, \\
&\int_{\mathbb{R}^{N}}\left|u\right|^{2} dx=c, \quad u \in H^{1}(\mathbb{R}^{N}),
\end{aligned}\right.
\end{equation}
which has attracted considerable attention. When $f(x,u)=|u|^{p-2}u$ (i.e., the limited problem of Eq.\eqref{initial}), the problem \eqref{yiban4} turns to
\begin{equation}\label{lim-initial}
\left\{
\begin{aligned}
&-\left(a+b\int_{\mathbb{R}^{N}}\left | \nabla u \right |^{2}dx\right)\Delta u+\lambda u=|u|^{p-2} u,   \quad \text { in } \mathbb{R}^{N}, \\
&\int_{\mathbb{R}^{N}} \left|u\right|^{2} dx=c, \quad u \in H^{1}(\mathbb{R}^{N}),
\end{aligned}
\right.
\end{equation}
where $a, b, c > 0$ are constants, $1\le N\le 3$, and $p\in \left(2,2^{*}\right)$. The energy functional of \eqref{lim-initial} is
\begin{equation}\label{Iwqu}
I_{\infty}(u)=\frac{a}{2}\int_{\mathbb{R}^{N}}\left|\nabla u\right|^{2} dx+\frac{b}{4} \left ( \int_{\mathbb{R}^{N}}\left|\nabla u\right|^{2} dx \right )^{2}-\frac{1}{p} \int_{\mathbb{R}^{N}}\left|u\right|^{p} dx.
\end{equation}
By the Galiardo-Nirenberg inequality \cite{r1} for any $p\in \left(2,2^{*}\right)$
\begin{eqnarray}\label{GN}
\|u\|_{p}\leq C_{N,p}\|\nabla u\|_{2}^{\gamma_{p}}\|u\|_{2}^{1-\gamma_{p}}
\end{eqnarray}
where $\gamma_{p}=\frac{N(p-2)}{2p},$ we can get $L^{2}$-critical exponent $\bar{p} =2+\frac{8}{N}$ of Kirchhoff problem. It is well known that Ye \cite{r6} obtained the sharp existence of global constraint minimizers for Eq.\eqref{lim-initial} in the case of $p\in (2, \bar{p})$. When $p\in (2+\frac{4}{N},\bar{p})$, Ye proved a local minimizer which is a critical point of $\left.I_{\infty}\right|_{S_{c}}$. By considering a global minimization problem
\begin{equation}\label{wqinf}
l_{\infty,c}:=\inf _{S_{c}}I_{\infty}(u),
\end{equation}
we have
\begin{equation}\label{jx0}
\begin{aligned}
\left\{
\begin{array}{l}
l_{\infty,c}\in (-\infty,0], \ \ if \ \ p\in (2,\bar{p}), \\
l_{\infty,c}=-\infty,\ \ \ \ \ \  if \ \ \ p\in (\bar{p} ,2^{*}), \\
\end{array}
\right.
\end{aligned}
\end{equation}
for any given $c>0$. We can see that the minimization method is not feasible for $p\in (\bar{p},2^{*})$. Then Ye proved the existence of normalized solutions by taking advantage of the Pohozaev constraint method in the case of $p\in (\bar{p} ,2^{*}).$ For the $L^{2}$-critical case of $\bar{p}=2+\frac{8}{N}$, Ye \cite{z14} showed the existence and mass concentration of critical points. Using some simple energy estimates instead of the concentration-compactness principles introduced in \cite{r6}, Zeng studied the existence and uniqueness of normalized solutions for $p\in \left (2,2^{*}\right)$ in \cite{z15}.

Additionally, Li, Luo and Yang \cite{L2} proved the existence and asymptotic properties of solutions to the following equation with combined nonlinearity
\begin{equation}\label{yiban5}
\left\{\begin{aligned}
&-\left(a+b\int_{\mathbb{R}^{N}}\left | \nabla u \right |^{2}dx\right)\Delta u+\lambda u=|u|^{p-2}u+\mu |u|^{q-2}u, \quad \text { in } \mathbb{R}^{3}, \\
&\int_{\mathbb{R}^{N}}\left|u\right|^{2} dx=c, \quad u \in H^{1}(\mathbb{R}^{N}),
\end{aligned}\right.
\end{equation}
where $a, b, c,\mu > 0$, $2<q<\frac{14}{3}< p \le 6$ or $\frac{14}{3}<q<p \le 6$, and showed a multiplicity result for the case of  $2 < q < \frac{10}{3}$ and $\frac{14}{3} < p < 6$, and the existence of ground state normalized solutions for $2 < q < \frac{10}{3} < p = 6$ or $\frac{14}{3}< q < p\le 6$. They also showed some asymptotic results on the obtained solutions.
For the case $\mu \le0$, in \cite{C6}, Carri$\tilde{a}$o, Miyagaki, and Vicente studied the ground states existence of Eq.\eqref{yiban5} for $2<q<2^{*}, p=2^{*}$ or $2<q \le \bar{p} <p<2^{*}$. For the nonautonomous problem, when $f(x,u)=|u|^{p-2}u+V(x)|u|^{q-2}u$, $N=3$, $p=\frac{14}{3},$ $q=4$ and $V \in L_{loc}^{\infty}(\mathbb{R}^{3})$, Ye \cite{z13} considered the existence of minimizers to the nonautonomous problem. Moreover, $V(x)$ satisfies
$$
V(x)\ge 0, \quad \lim_{\left | x \right |  \to \infty }V(x)=0.
$$
By the concentration compactness principle, if $b<b_{0}$, Ye showed that there exists $a_{0}, c_{0}>0$ such that the above problem has a minimizer for all $a <a_{0}$ and $c<c_{0}$. Additionally, when $f(x,u)=K(x)f(u)$, Chen and Tang \cite{z16} considered the existence of ground state solutions, where $K(x) \in C(\mathbb{R}^{3},\mathbb{R}^{+})$ and $f(u)$ is $L^{2}$-supercritical. Other results about normalized solutions of Kirchhoff equation in a more general form can refer to \cite{z16,z17,C7,Z1}.

Motivated by the results above, when $\mu$ of Eq.\eqref{yiban5} is replaced by a potential function $h(x)$ and $1\le q< 2$, there is no results in studying normalized solutions of such nonautonomous Kirchhoff equations with a small perturbation. In the present paper, we first obtain the normalized solution of this type equation, which can be seen as extension of some known results in the literature.

Let us now outline the main strategy to prove the three results of this paper under different assumptions on $h(x)$. Firstly, we treat the mass-subcritical case $2<p<2+\frac{4}{N}$: for any $c>0$, we set
\begin{equation}\label{inf}
l_{c}:=\inf _{S_{c}}I(u).
\end{equation}
It is standard that the minimizers of $l_{c}$ are critical points of $\left.I\right|_{S_{c}}$. We introduce the following assumptions on $h(x)$.
$$
\left(\mathbf{h}_{\mathbf{1}}\right)\ \ h \in L^{\frac{2}{2-q}}(\mathbb{R}^{N}) \text { and }h(x)>0\ \  \text {on a set with positive measure.}
$$
Now we state the main results of this paper:\\

\textbf{Theorem 1.1.}\ \ Suppose $1\le N\le3$, $2<p<2+\frac{4}{N}$ and $h(x)\ge 0$ satisfy $\left(\mathbf{h}_{\mathbf{1}}\right)$. Then for all $c>0$, $l_{c}$ has a minimizer, hence Eq.\eqref{initial} has a normalized ground state solution.\\

\textbf{Remark 1.1.}\ \ Notice that the minimizer obtained in Theorem 1.1 is a global minimizer rather than a local minimizer. It is easy to find that the energy functional is coercive on $S_{c}$, which hints each minimizing sequence $\{u_{n}\}$ is bounded on $S_{c}$. The main difficulty of proof is to show the minimizing sequence $\{u_{n}\}$ converge strongly to $u\ne 0$ in $H^{1}(\mathbb{R}^{N})$. The key step is to establish the inequality $l_{c_{1}+c_{2}}\leq l_{c_{1}}+l_{\infty,c_{2}}$ for $c_{1}, c_{2}>0$ (see Lemma 2.2), which is crucial to recover the compactness.\\

Next, while addressing the $L^{2}$-supercritical case, the functional is unbounded from below on $S_{c}$, thus the minimizing approach on $S_{c}$ is not valid anymore. Ye \cite{r6} proved that $l_{\infty,c}=-\infty$ for all $c>0$ if $p\in \left(2+\frac{8}{N},2^{*}\right)$, and Ye proved the existence of one normalized solution by a suitable submanifold of $S_{c}$. In this paper, after the appearance of a very small perturbation term, we want to show that the energy functional $I$ has a mountain pass geometry and show the existence of a mountain pass solution with positive energy level for $p\in \left(2+\frac{8}{N},2^{*}\right)$. We require the perturbation $h(x)$ to have a higher regularity. We need to assume that:
$$
\left(\mathbf{h}_{\mathbf{2}}\right)\ \ h \in L^{\frac{p}{p-q}}(\mathbb{R}^{N})\cap C^{1}(\mathbb{R}^{N}), \ \ \langle\nabla h, x\rangle  \in L^{\frac{2}{2-q} }(\mathbb{R}^{N})\ \text{and}\ h(x)\ge0.
$$
We have the following result.\\

\textbf{Theorem 1.2.}\ \ Suppose $1\le N\le3$, $2+\frac{8}{N}<p<2^{*}$ and $h(x)$ satisfy $\left(\mathbf{h}_{\mathbf{2}}\right)$. Let $c>0$ be fixed. 
Moreover,
\begin{equation}\label{Lcondi1}
\left \| h \right \|_{\frac{p}{p-q}}<\frac{aq\left(p \gamma_{p}-2\right)}{2\mathcal{C}_{_{N, p}}^{q}  \gamma_{p}(p-q)}\left(\frac{ap\left(2-q\gamma_{p}\right)}{2\gamma_{p}(p-q) \mathcal{C}_{N, p}^{p}}\right)^{\frac{2-q\gamma_{p}}{p\gamma_{p}-2}} c^{-\frac{(1-\gamma_{p})(p-q)}{p\gamma_{p}-2}},
\end{equation}
\begin{equation}\label{Lcondi2}
\|\nabla h \cdot x\|_{\frac{2}{2-q}}<\frac{q(2p-Np+2N)}{p-2}m_{c}c^{-\frac{q}{2}}.
\end{equation}
Then Eq.\eqref{initial} has a mountain pass solution $u$ at a positive energy level.\\

\textbf{Remark 1.2.}\ \ We are going to use the minimax characterization to find a critical point. Although the mountain-pass geometry of the functional on $S_{c}$ can be obtained easily, unfortunately the boundedness of the obtained $(PS)$ sequence is not yet clear. In this paper, we adopt a similar idea in \cite{r4} and construct an auxiliary map $\tilde{I}(t, u):=I(t\star u),$ which on $\mathbb{R}\times S_{c}$ has the same type of geometric structure as $I$ on $S_{c}$. Besides, the $(PS)$ sequence of $I$ satisfies the additional condition(see Lemma 3.5), which is the key ingredient to obtain the boundedness of the $(PS)$ sequence.\\

Finally, we will discuss $h(x)\le0$, the problem becomes more delicate and difficult. Although the mountain pass structure by Jeanjean \cite{r4} is destroyed, Bartsch et al.\cite{r18} established a new variational principle exploiting the Pohozaev identity. For convenience's sake, we define $\bar{h}(x):=-h(x)\ge 0$. Next, we state our basic assumptions on $\bar{h}(x)$.\\

$\left(\mathbf{h}_{\mathbf{3}}\right)\ \bar{h}(x)\in L^{\frac{2}{2-q} }(\mathbb{R}^{N})\cap C^{1}(\mathbb{R}^{N}),\ \langle\nabla \bar{h}(x), x\rangle \in L^{\frac{2}{2-q} }(\mathbb{R}^{N})$ and $\bar{h}(x)\ge0$. For some constants $\Upsilon> 0$, $\bar{h}(x)$ satisfies
$$\left|x\cdot\nabla \bar{h}(x) \right |\le \Upsilon \bar{h}(x).$$ \\

\textbf{Theorem 1.3.}\ \ Assume $1\le N\le3$, $2+\frac{8}{N}<p<2^{*}$. If $\left(\mathbf{h}_{\mathbf{3}}\right)$ holds and $\bar{h}(x)$ satisfies
\begin{equation}\label{th13}
0< \|\bar{h}\|_{\frac{2}{2-q}}< \min \left \{ 1,\frac{2p\left(1-\gamma_{p}\right)}{2(p-q)+(p-2)\Upsilon}\right \}\cdot  \frac{qm_{c}}{c^{\frac{q}{2}}}.
\end{equation}
Then Eq.\eqref{initial} has a couple of solution $\left(u_,\lambda\right) \in H^{1}(\mathbb{R}^{N})\times \mathbb{R}$ and $\lambda>0$.\\

\textbf{Remark 1.3.}\ \ Indeed, when $h(x)\leq0$, the problem is made more difficult by the simultaneous appearance of negative potential and nonlocal term. We refer to Bartsch et al. in \cite{r18} constructing a suitable linking geometry method to obtain the existence of bound state solutions with high Morse index. The crucial step is to estimate the minimax level $m_{c}<L_{h,c}<2m_{c}$ (see Lemma 4.3 and Lemma 4.5) to recover the compactness.\\

\textbf{Notations:}\ \ We introduce some notations that will clarify what follows:\\
$\bullet\ \ H^{1}(\mathbb{R}^{N})$ is the usual Sobolev space with the norm $\left\|u\right\|=\left(\int_{\mathbb{R}^{N}}\left|\nabla u\right|^{2}+\left | u \right|^{2} dx \right)^{\frac{1}{2}}$.\\
$\bullet\ \ L^{p}(\mathbb{R}^{N})$ with $p \in [1,\infty)$ is the Lebesgue space with the norm $\left \|u\right \|_{p}=\left ( \int_{\mathbb{R}^{N}}\left| u\right|^{p}dx \right)^{\frac{1}{p} }$.\\
$\bullet$\ \ The arrows $'\rightharpoonup'$ and $'\to'$ denote the weak convergence and strong convergence, respectively.\\
$\bullet\ \ C,C_{i}$ denote positive constants which may vary from line to
line.\\
$\bullet\ \ (t\star u)(x):=t^{\frac{N}{2}}u(tx)$ for $t\in \mathbb{R}^{+}$ and $u\in H^{1}(\mathbb{R}^{N})$.

\section{Proof of Theorem 1.1}

\ \ \ \ \ \ In this section, for $2<p<2+\frac{4}{N}$ and $h(x)\ge0$ we prove Theorem 1.1. By the Gagliardo-Nirenberg inequality \eqref{GN},
\begin{equation}\label{0}
\begin{aligned}
I(u)&=\frac{a}{2}\|\nabla u\|_{2}^{2}+\frac{b}{4}\|\nabla u\|_{2}^{4}-\frac{1}{p}\|u\|_{p}^{p}-\frac{1}{q}\int_{\mathbb{R}^{N}}h(x)\left|u \right|^{q} dx\\
 &\geq \frac{a}{2}\|\nabla u\|_{2}^{2}+\frac{b}{4}\|\nabla u\|_{2}^{4}-\frac{1}{p} \mathcal{C}_{N, p}^{p}\|\nabla u\|_{2}^{p\gamma_{p}}\|u\|_{2}^{p\left(1-\gamma_{p}\right)}- \frac{1}{q}\|h\|_\frac{2}{2-q}\|u\|_{2}^{q},\\
\end{aligned}
\end{equation}
thus $I$ is bounded from below on $S_{c}$ since $0<p\gamma_{p}<2$.

For $1\le N\le3$ and $2<p<2+\frac{4}{N}$, the existence and uniqueness of positive normalized solutions of the limited problem \eqref{lim-initial} has been studied  in \cite{r6}. In order to find the minimizer of $I$ on $S_{c}$, firstly we state some fundamental properties of $l_{\infty,c}$, which will be crucial to recover the compactness later on. The proof of next lemma can be referred to Theorem 1.1 and Lemma 2.5 in \cite{C7}.\\



\textbf{Lemma 2.1.}\ \ Suppose $1\le N\le3$ and $2<p<2+\frac{4}{N}$. Then for all $c>0$, we have\\
(i)the strict sub-additivity for $l_{\infty,c}$, i.e.,
$$
l_{\infty,c_{1}+c_{2}}<l_{\infty,c_{1}}+l_{\infty,c_{2}} \ \ \text{for} \ \ c_{1}, c_{2}>0;
$$
(ii)the limited problem \eqref{lim-initial} has a couple of ground state solution $(u_{\infty },\lambda _{c})\in H^{1}(\mathbb{R}^{N}) \times \mathbb{R}$, i.e.,
$$
l_{\infty,c}=\inf_{S_{c}}I_{\infty}(u)=I_{\infty}(u_{\infty })<0.
$$


Next, we introduce the inequality $l_{c_{1}+c_{2}}\leq l_{c_{1}}+l_{\infty,c_{2}}$ which plays a crucial role in proving the convergence of the minimizing sequence.\\

\textbf{Lemma 2.2.}\ \ Suppose $2<p<2+\frac{4}{N}$ and  $h(x)$ satisfy $\left(\mathbf{h_{1}}\right)$, then the following holds\\
(i)$-\infty<l_{c}<l_{\infty,c}<0$ for $c>0$;\\
(ii)$l_{c_{1}+c_{2}}\leq l_{c_{1}}+l_{\infty,c_{2}}$ for $c_{1}, c_{2}>0$.

\textbf{Proof.}\ \ (i) It is obvious that $l_{c}>-\infty$ by \eqref{0}. Moreover, by Lemma 2.1, we have
$$
\begin{aligned}
l_{c}&\leq I\left(u_{\infty}\right) \\
&=\frac{a}{2} \int_{\mathbb{R}^{N}}\left|\nabla u_{\infty}\right|^{2}dx+\frac{b}{4} \left ( \int_{\mathbb{R}^{N}}\left|\nabla u_{\infty}\right|^{2} dx \right )^{2}-\frac{1}{p} \int_{\mathbb{R}^{N}}\left|u_{\infty}\right|^{p} dx-\frac{1}{q} \int_{\mathbb{R}^{N}} h\left | u_{\infty} \right |^{q} dx \\
& <I_{\infty}\left(u_{\infty}\right) \\
& =l_{\infty,c}<0,
\end{aligned}
$$
since $u_{\infty}>0$ and $h(x)$ satisfies $\left(\mathbf{h_{1}}\right)$.

(ii) For any $\varepsilon>0$, $c=c_{1}+c_{2}$, we can find  $\varphi_{\varepsilon}, \psi_{\varepsilon} \in \mathcal{C}_{0}^{\infty}(\mathbb{R}^{N})$ such that
$$
\begin{aligned}
\varphi_{\varepsilon} \in S_{c_{1}} , \quad I\left(\varphi_{\varepsilon}\right)<l_{c_{1}} +\frac{\varepsilon}{2}, \\
\psi_{\varepsilon} \in S_{c_{2} }, \quad I_{\infty } \left(\psi_{\varepsilon}\right)<l_{\infty ,c_{2}}+\frac{\varepsilon}{2}.
\end{aligned}
$$
Let $u_{\varepsilon,n}(x):=\varphi_{\varepsilon}(x)+\psi_{\varepsilon}\left(x-n \boldsymbol{e}_{1}\right)$, where $\boldsymbol{e}_{1}$ is the unit vector $(1,0, \cdots)$ in $\mathbb{R}^{N}$. Since $\varphi_{\varepsilon}$ and $\psi_{\varepsilon}$ have compact support, we see that $u_{\varepsilon,n} \in S_{c}$ for large $n$ and that
$$
l_{c}\leq I\left(u_{\varepsilon, n}\right)=I\left(\varphi_{\varepsilon}\right)+I\left(\psi_{\varepsilon}\left(x-n \boldsymbol{e}_{1}\right)\right) .
$$
Moreover, thanks to $h \in L^{\frac{2}{2-q}}(\mathbb{R}^{N})$, we have that $\int_{\mathbb{R}^{N}}h(x)\psi_{\varepsilon}^{q} \left(x-n \boldsymbol{e}_{1}\right)dx\to 0$ as $n\to \infty$, hence $I(\psi_{\varepsilon}\left(\cdot -n \boldsymbol{e}_{1}\right))\to I_{\infty } (\psi_{\varepsilon})$ as $n\to \infty$.
It follows that
$$
\begin{aligned}
 l_{c}& \leq \limsup _{n \rightarrow \infty} I\left(u_{\varepsilon, n}\right) \\
& =\limsup _{n \rightarrow \infty}\left(I\left(\varphi_{\varepsilon}\right)+I\left(\psi_{\varepsilon}\left(\cdot-n \boldsymbol{e}_{1}\right)\right)\right) \\
&=I\left(\varphi_{\varepsilon}\right)+I_{\infty}\left(\psi_{\varepsilon}\right) \\
&<l_{c_{1}}+l_{\infty,c_{2}}+\varepsilon .
\end{aligned}
$$
Passing to the limit, thus $l_{c}\leq l_{c_{1}}+l_{\infty,c_{2}}$ since $\varepsilon>0$ is arbitrary.\ \ \ $\Box$\\

Let $\{u_{n}\}\subset S_{c}$ be a minimizing sequence for $l_{c}$. By \eqref{0}, we know that $I(u)$ is coercive on $S_{c}$ and deduce that $\{u_{n}\}$ is bounded in $H^{1}(\mathbb{R}^{N})$. Thus, there exists a subsequence such that $u_{n}\rightharpoonup u_{0}$ and
$$I(u_{0})\le \liminf _{n \rightarrow \infty}I(u_{n})= l_{c},\quad  c_{1}:=\left \| u _{0} \right \|_{2}^{2}\le \left \| u _{n} \right \|_{2}^{2}=c.$$
We need to prove $I(u_{0})=l_{c}$ and $\left \| u _{0} \right \|_{2}^{2}=c$. Now we argue by contradiction to prove this.\\

\textbf{Lemma 2.3.}\ \ Suppose $2<p<2+\frac{4}{N}$ and $h(x)$ satisfying $\left(\mathbf{h_{1}}\right)$. Then every minimizing sequence for $l_{c}$ has a strong convergent subsequence in $L^{2}(\mathbb{R}^{N})$.

\textbf{Proof.}\ \ We argue by contradiction and assume that $c_{1}<c$. We divide the proof into four steps.

\textbf{Step1}: there exists $\{y_{n}\}\subset \mathbb{R}^{N}$ and $\mu_{0} \in H^{1}(\mathbb{R}^{N}) \backslash\{0\}$ such that
\begin{equation}\label{x1}
\left|y_{n}\right| \rightarrow \infty, \ \ u_{n}\left(\cdot+y_{n}\right)\rightharpoonup \mu_{0} \quad \text { in } H^{1}(\mathbb{R}^{N}).
\end{equation}
First, we show by contradiction that
\begin{equation}\label{x5}
\delta_{0}:=\liminf _{n \rightarrow \infty} \sup _{y \in \mathbb{R}^{N}} \int_{B_{1}(y)}\left|u_{n}-u_{0}\right|^{2} dx>0,
\end{equation}
where $B_{1}(y)=\{x \in \mathbb{R}^{N}: |x-y|\leq 1\}$. Suppose on the contrary that $\delta_{0}=0$. Then, $u_{n} \rightarrow u_{0}$ strongly in $L^{p}(\mathbb{R}^{N})$ (see \cite{r3}). Since $u_{n} \rightharpoonup u_{0}$ in $H^{1}(\mathbb{R}^{N}),  h \in L^{\frac{2}{2-q}}(\mathbb{R}^{N})$, we see that $\int_{\mathbb{R}^{N}} h \left|u_{n} \right |^{q} dx \rightarrow \int_{\mathbb{R}^{N}} h \left|u_{0} \right |^{q} dx$. Combined with Lemma 2.1 (ii), for $c-c_{1}>0$, we have that
$$
\begin{aligned}
l_{c} & =I\left(u_{n}\right)+o(1) \\
& =I\left(u_{0}\right)+I\left(u_{n}-u_{0}\right)+o(1) \\
& =I\left(u_{0}\right)+\frac{a}{2} \int_{\mathbb{R}^{N}}\left|\nabla\left(u_{n}-u_{0}\right)\right|^{2} dx+\frac{b}{4} \left ( \int_{\mathbb{R}^{N}}\left|\nabla\left(u_{n}-u_{0}\right)\right|^{2} dx \right )^{2}  +o(1) \\
& >l_{c_{1}}+l_{\infty,c-c_{1}},
\end{aligned}
$$
which is a contradiction with Lemma 2.2 (ii). Therefore, \eqref{x5} holds.
From \eqref{x5} and $u_{n} \rightarrow u_{0}$ in $L_{loc}^{2}(\mathbb{R}^{N})$, we can find $\{y_{n}\} \subset \mathbb{R}^{N}$ such that  $\int_{B_{1}(y_{n})}\left|u_{n}-u_{0}\right|^{2} dx \rightarrow c_{0}>0$ and $\left|y_{n}\right| \rightarrow \infty$. Let $u_{n}\left(\cdot+y_{n}\right) \rightharpoonup \mu_{0}$ weakly in $H^{1}(\mathbb{R}^{N})$. Note that $\mu_{0}\neq 0$ since $c_{0}>0$. Therefore, $\{y_{n}\}$ and $\mu_{0}$ satisfy \eqref{x1}. Thus, the proof of Step 1 is complete.

\textbf{Step2}: We show that $\{y_{n}\}$ and  $\left(u_{0}, \mu_{0}\right)$ satisfy
\begin{equation}\label{x2}
\lim _{n \rightarrow \infty}\|u_{n}-u_{0}-\mu_{0}\left(\cdot-y_{n}\right)\|_{2}^{2}=0.
\end{equation}
Since $\left|y_{n}\right| \rightarrow \infty$, we have that
\begin{equation}\label{s1}
\begin{aligned}
\|u_{n}-u_{0}-\mu_{0}(\cdot-y_{n})\|_{2}^{2}= & \|u_{n}\|_{2}^{2}+\|u_{0}\|_{2}^{2}+\|\mu_{0}\|_{2}^{2} \\
& -2\left\langle u_{n}, u_{0}\right\rangle_{L^{2}}-2\left\langle u_{n}\left(\cdot+y_{n}\right), \mu_{0}\right\rangle_{L^{2}}+o(1) \\
= & \|u_{n}\|_{2}^{2}-\|u_{0}\|_{2}^{2}-\|\mu_{0}\|_{2}^{2}+o(1).
\end{aligned}
\end{equation}
According to \eqref{s1}, we could let $\delta_{1}:=\lim _{n \rightarrow \infty}\|u_{n}-u_{0}-\mu_{0}(\cdot-y_{n})\|_{2}^{2}$. Then, we have $\delta_{1}=c-c_{1}-c_{2}$, where $c_{2}:= \|\mu_{0}\|_{2}^{2}$. We want to show that $\delta_{1}=0$. Suppose on the contrary that  $\delta_{1}>0$, by direct calculations we have
\begin{equation}\label{s2}
\begin{aligned}
&\|\nabla u_{n}\|_{2}^{2}-\|\nabla u_{0}\|_{2}^{2}-\|\nabla \mu_{0}(\cdot-y_{n})\|_{2}^{2}-\|\nabla(u_{n}-u_{0}-\mu_{0}(\cdot-y_{n}))\|_{2}^{2} \\
\quad&=-2\|\nabla u_{0}\|_{2}^{2}-2\|\nabla \mu_{0}\|_{2}^{2}+2\left\langle\nabla u_{n}, \nabla u_{0}\right\rangle_{L^{2}}+2\left\langle\nabla u_{n}\left(\cdot+y_{n}\right), \nabla\mu_{0}\right\rangle_{L^{2}} \\
\quad&=o(1).
\end{aligned}
\end{equation}
From the Brezis-Lieb Lemma, we have
\begin{equation}\label{s3}
\begin{aligned}
\int_{\mathbb{R}^{N}}\left|u_{n}\right|^{p} dx= & \int_{\mathbb{R}^{N}}\left|u_{0}\right|^{p}dx+\int_{\mathbb{R}^{N}}\left|\mu_{0}\left(\cdot-y_{n}\right)\right|^{p} dx \\
& +\int_{\mathbb{R}^{N}}\left|u_{n}-u_{0}-\mu_{0}\left(\cdot-y_{n}\right)\right|^{p} dx+o(1) .
\end{aligned}
\end{equation}
Similarly,
\begin{equation}\label{s4}
\begin{aligned}
\int_{\mathbb{R}^{N}}h\left | u_{n} \right |^{q}  \mathrm{~d} x=&\int_{\mathbb{R}^{N}}h\left | u_{0} \right |^{q} dx+\int_{\mathbb{R}^{N}}h\left | \mu_{0}\left(\cdot-y_{n}\right)\right|^{q} dx\\
&+\int_{\mathbb{R}^{N}} h\left | \left(u_{n}-u_{0}-\mu_{0}\left(\cdot-y_{n}\right)\right) \right|^{q} dx+o(1).
\end{aligned}
\end{equation}
Combining \eqref{s2}-\eqref{s4}, we have
\begin{equation}\label{s5}
I\left(u_{n}\right)-I\left(u_{0}\right)-I\left(\mu_{0}\left(\cdot-y_{n}\right)\right)-I\left(u_{n}-u_{0}-\mu_{0}\left(\cdot-y_{n}\right)\right)=o(1).
\end{equation}
Since $u_{n} \rightharpoonup u_{0}$ in $H^{1}(\mathbb{R}^{N})$, $\left|y_{n}\right| \rightarrow \infty$ and $h \in L^{\frac{2}{2-q}}\left(\mathbb{R}^{N}\right)$, we have
\begin{equation}\label{s6}
\int_{\mathbb{R}^{N}} h\left | u_{n}-u_{0}-\mu_{0}\left(\cdot-y_{n}\right) \right |^{q} dx \rightarrow 0.
\end{equation}
Recalling that $l_{\infty,c}$ is continuous with respect to $c> 0$ (see \cite{r16}, Theorem 2.1), we have that
\begin{equation}\label{s7}
\begin{aligned}
\begin{array}{l}
\liminf _{n \rightarrow \infty} I\left(u_{n}-u_{0}-\mu_{0}\left(\cdot-y_{n}\right)\right) \\
\quad=\liminf _{n \rightarrow \infty} I_{\infty}\left(u_{n}-u_{0}-\mu_{0}\left(\cdot-y_{n}\right)\right) \\
\quad \geq l_{\infty,\delta_{1}},
\end{array}
\end{aligned}
\end{equation}
and
\begin{equation}\label{s8}
\liminf _{n \rightarrow \infty} I\left(\mu_{0}\left(\cdot-y_{n}\right)\right) \geq l_{\infty,c_{2}}.
\end{equation}
Hence by \eqref{s5}-\eqref{s8}, we have
\begin{equation}\label{s9}
l_{c} \geq l_{c_{1}}+l_{\infty,c_{2}}+l_{\infty,\delta_{1}}.
\end{equation}
However, using Lemma 2.1 (i), for any $c_{2},\delta_{1}>0$, there exist $l_{\infty,c_{2}+\delta_{1}}< l_{\infty,c_{2}}+l_{\infty,\delta_{1}}$. Hence we also have
\begin{equation}\label{s10}
\begin{aligned}
l_{c} &\geq l_{c_{1}}+l_{\infty,c_{2}}+l_{\infty,\delta_{1}} \\
&>l_{c_{1}}+l_{\infty,c_{2}+\delta_{1}}\\
&\geq l_{c_{1}+c_{2}+\delta_{1}} \\
& =l_{c} .
\end{aligned}
\end{equation}
This gives a contradiction and thus we have that $\delta_{1}=0$.

\textbf{Step3}: Moreover, the following hold
\begin{equation}\label {x3}
I \left(u_{0}\right)=l_{c_{1}}, \ \ I_{\infty}\left(\mu_{0}\right)=l_{\infty,c_{2}},
\end{equation}
and
\begin{equation}\label {x4}
l_{c}=l_{c_{1}}+l_{\infty,c_{2}}.
\end{equation}
By \eqref{s5}-\eqref{s8} and $\delta_{1}=0$, we have that
\begin{equation}\label{s11}
\begin{aligned}
l_{c} & =\lim _{n \rightarrow \infty} I\left(u_{n}\right) \\
& =\liminf _{n \rightarrow \infty}\left(I\left(u_{0}\right)+I\left(\mu_{0}\left(\cdot+y_{n}\right)\right)\right) \\
& \geq I\left(u_{0}\right)+I_{\infty}\left(\mu_{0}\right) \\
& \geq l_{c_{1}}+l_{\infty,c_{2}}.
\end{aligned}
\end{equation}
Combined with Lemma 2.2 (ii), we see that  $l_{c}=l_{c_{1}}+l_{\infty,c_{2}}$. $I\left(u_{0}\right)=l_{c_{1}}$ and $I_{\infty}\left(\mu_{0}\right)=l_{\infty,c_{2}}$. Thus, Step 3 is proved.

\textbf{Step4}: Now, we prove the precompactness of minimizing sequence, i.e., $u _{n}\to u _{0}$ in $L^{2}(\mathbb{R}^{N})$.
 \\
We can suppose that $\{u _{n}\}$ are nonnegative. Using the strong maximum principle, we have $u_{0}, \mu_{0}>0$ and $h(x)>0$ on a set with positive measure, we have that
$$
\int_{\mathbb{R}^{N}} h\left |\sqrt{u_{0}^{2}+\mu_{0}^{2}} \right |^{q}\mathrm{~d} x>\int_{\mathbb{R}^{N}} h\left | u_{0}\right |^{q}\mathrm{~d}x.
$$
Combine  with the two following inequalities:
\begin{equation}\label{qq5}
\int_{\mathbb{R}^{N}}\left|\nabla \sqrt{u_{0}^{2}+\mu_{0}^{2}}\right|^{2} dx \leq \int_{\mathbb{R}^{N}}(\left|\nabla u_{0}\right|^{2}+\left|\nabla \mu_{0}\right|^{2}) dx,
\end{equation}
\begin{equation}\label{q5}
\int_{\mathbb{R}^{N}}\left|\sqrt{u_{0}^{2}+\mu_{0}^{2}}\right|^{p} dx \geq \int_{\mathbb{R}^{N}}\left(\left|u_{0}\right|^{p}+\left|\mu_{0}\right|^{p}\right) dx.
\end{equation}
So we have
\begin{equation}\label{q6}
\begin{aligned}
l_{c} & \leq I\left(\sqrt{u_{0}^{2}+\mu_{0}^{2}}\right)\\
& =\frac{a}{2} \int_{\mathbb{R}^{N}}\left|\nabla \sqrt{u_{0}^{2}+\mu_{0}^{2}}\right|^{2} dx+\frac{b}{4}\left(\int_{\mathbb{R}^{N}}\left|\nabla \sqrt{u_{0}^{2}+\mu_{0}^{2}}\right|^{2} dx\right )^{2}\\
&-\frac{1}{p} \int_{\mathbb{R}^{N}}\left|\sqrt{u_{0}^{2}+\mu_{0}^{2}}\right|^{p} dx-\frac{1}{q}\int_{\mathbb{R}^{N}}h \left | \sqrt{u_{0}^{2}+\mu_{0}^{2}} \right |^{q} dx \\
& <I\left(u_{0}\right)+I_{\infty}\left(\mu_{0}\right) \\
& =l_{c_{1}}+l_{\infty,c-c_{1}} \\
& =l_{c},
\end{aligned}
\end{equation}
which is a contradiction. Thus the proof of Lemma 2.3 is completed.\ \ \ $\Box$\\

\textbf{Proof of Theorem 1.1.}\ \ From Lemma 2.3, the minimizing sequence $\{u_{n}\}$ satisfy $u_{n}\rightarrow u_{0}$ in  $L^{2}\left(\mathbb{R}^{N}\right)$ and $l_{c}=I\left(u_{0}\right)$, $c=c_{1}$. Since $\{u_{n}\}\subset S_{c}$ is the minimizing sequence of $l_{c}$, we have $\left.d I\right|_{S_{c}}\left(u_{n}\right) \rightarrow 0$ and there exists a sequences of real numbers $\{\lambda_{n}\}$ such that
\begin{equation}\label{s12}
I^{\prime}\left(u_{n}\right)[\varphi]+\lambda_{n} \int_{\mathbb{R}^{N}} u_{n} \varphi dx \rightarrow 0, \quad as \ \ n \rightarrow \infty,
\end{equation}
for every $\varphi \in H^{1}(\mathbb{R}^{N})$. Hence by \eqref{s12}, we have that
\begin{equation}\label{q7}
\left\{\begin{aligned}
&-\left ( a+b\int _{\mathbb{R}^{N}}\left | \nabla u_{0} \right |^{2}dx \right ) \Delta u_{0}+\bar{\lambda} u_{0}=|u_{0}|^{p-2} u_{0}+h(x)\left | u_{0} \right |^{q-2}u_{0}  \quad \text { in } \mathbb{R}^{N}, \\
&\int_{\mathbb{R}^{N}} \left|u_{0}\right|^{2} dx=c.
\end{aligned}\right.
\end{equation}
Notice that $h(x)\geq 0$, then by the maximum principle,  $u_{0}>0$ and we finish the proof of Theorem 1.1.\ \ \ $\Box$\\

\section{Proof of Theorem 1.2}

\ \ \ \ \ \ In this section, we study the mass-supercritical and Sobolev-subcritical case: $2+\frac{8}{N} <p< 2^{*}$, $1\le N \le 3$, $h(x)\ge0$. Firstly, we show that the energy functional $I$ possesses a mountain pass geometry, which implies the existence of the $(PS)$ sequence. Next, we prove that the limit of the sequence of the Lagrange multipliers related to the $(PS)$ sequence is positive. Then by applying the splitting lemma, we recover the compactness for this sequence, which yields the existence of solutions for Eq.\eqref{initial}.

In order to study the behavior of $(PS)$ sequence, we introduce the splitting lemma which plays a crucial role in overcoming the lack of compactness. For $\lambda>0$ we set
$$
I_{\lambda}(u)=\frac{a}{2} \int_{\mathbb{R}^{N}}|\nabla u|^{2} d x+\frac{b}{4}\left(\int_{\mathbb{R}^{N}}|\nabla u|^{2} dx\right)^{2}+\frac{1}{2} \int_{\mathbb{R}^{N}}\lambda u^{2} dx-\frac{1}{p} \int_{\mathbb{R}^{N}}|u|^{p} dx-\frac{1}{q} \int_{\mathbb{R}^{N}} h\left | u \right |^{q} dx
$$
and
$$
I_{\infty,\lambda}(u)=\frac{a}{2} \int_{\mathbb{R}^{N}}|\nabla u|^{2} dx+\frac{b}{4}\left(\int_{\mathbb{R}^{N}}|\nabla u|^{2}dx\right)^{2}+\frac{1}{2} \int_{\mathbb{R}^{N}}\lambda u^{2} dx-\frac{1}{p} \int_{\mathbb{R}^{N}}|u|^{p}dx.
$$

\textbf{Lemma 3.1.}\ \ Let $\{u_{n}\}\subset H^{1}(\mathbb{R}^{N})$ be a $(PS)$ sequence for $I_{\lambda}$ such that $u_{n}\rightharpoonup u$ in $H^{1}(\mathbb{R}^{N})$ and $\lim_{n\to \infty}\left \|\nabla u_{n}\right \|_{2}^{2}=A^{2}$. Then there exists an integer $k\geq 0$, $k$ nontrivial solutions  $\omega^{1}, \cdots, \omega^{k} \in H^{1}(\mathbb{R}^{N})$ to the following problem
\begin{equation}\label{spit1}
-\left(a+bA^{2}\right) \Delta \omega+\lambda \omega=|\omega|^{p-2}\omega,
\end{equation}
and $k$ sequences $\{y_{n}^{j}\} \subset \mathbb{R}^{N}, 1\leq j \leq k$, such that as $n \rightarrow \infty,|y_{n}^{j}| \rightarrow \infty, |y_{n}^{j_{1}}-y_{n}^{j_{2}}| \rightarrow \infty$ for each $1 \leq j_{1}, j_{2} \leq k, j_{1} \neq j_{2}$, and
\begin{equation}\label{spit2}
\left\|u_{n}-u-\sum_{j=1}^{k} \omega^{j}\left(\cdot-y_{n}^{j}\right)\right\| \rightarrow 0,
\end{equation}
\begin{equation}\label{spit3}
A^{2}=\|\nabla u\|_{2}^{2}+\sum_{j=1}^{k}\left\|\nabla \omega^{j}\right\|_{2}^{2},
\end{equation}
\begin{equation}\label{2fs}
\left\|u_{n}\right\|_{2}^{2}=\|u\|_{2}^{2}+\sum_{j=1}^{k}\left\|w^{j}\right\|_{2}^{2}+o(1),
\end{equation}
and
\begin{equation}\label{spit4}
I_{\lambda}\left(u_{n}\right) \rightarrow J_{h,\lambda}(u)+\sum_{j=1}^{k} J_{\infty,\lambda}\left(\omega^{j}\right),
\end{equation}
as $n\to \infty$ where
$$
\begin{aligned}
J_{h,\lambda}(u):= &\left(\frac{a}{2}+\frac{b A^{2}}{4}\right)\int_{\mathbb{R}^{N}}|\nabla u|^{2} dx+\frac{\lambda}{2} \int_{\mathbb{R}^{N}} u^{2} dx\\
& -\frac{1}{p} \int_{\mathbb{R}^{N}}|u|^{p} dx-\frac{1}{q} \int_{\mathbb{R}^{N}} h\left | u \right |^{q} dx
\end{aligned}
$$
and
$$
J_{\infty,\lambda}(u):=\left(\frac{a}{2}+\frac{b A^{2}}{4}\right) \int_{\mathbb{R}^{N}}|\nabla u|^{2} dx+\frac{\lambda}{2} \int_{\mathbb{R}^{N}} u^{2} dx-\frac{1}{p} \int_{\mathbb{R}^{N}}|u|^{p}dx.
$$\\

\textbf{Lemma 3.2.}\ \ Let $X$ be a Hilbert manifold and let $F \in C^{1}(X,\mathbb{R})$  be a given functional. Let $K \subseteq X$ be compact and consider a subset.
$$
\mathcal{E} \subset\{E \subset X: E \text { is compact, } K \subset E\},
$$
which is invariant with respect to deformations leaving  $K$ fixed. Assume that
$$
\max _{u\in K}F(u)<c:=\inf _{E \in \mathcal{E}} \max _{u\in E}F(u) \in \mathbb{R}.
$$
Let  $\sigma_{n} \in \mathbb{R} $  be such that  $\sigma_{n} \rightarrow 0$  and  $E_{n} \in \mathcal{E}$  be a sequence such that
$$
c \leq \max _{u \in E_{n}}F(u)<c+\sigma_{n}.
$$
Then there exists a sequence $v_{n}\in X$ such that\\
1.  $c \leq F\left(v_{n}\right)<c+\sigma_{n} ,$\\
2.  $\left\|\nabla_{X} F\left(v_{n}\right)\right\|<\tilde{c} \sqrt{\sigma_{n}} ,$\\
3.  $\operatorname{dist}\left(v_{n}, E_{n}\right)<\tilde{c} \sqrt{\sigma_{n}} ,$\\
for some constant $\tilde{c}>0$.\\

We shall that $I$ on $S_{c}$ possesses a kind of mountain-pass geometrical structure. To this aim, we establish two preliminary lemmas.\\

\textbf{Lemma 3.3.}\ \ Assume that $h \in L^{\frac{p}{p-q}}(\mathbb{R}^{N})$ and let $u \in S_{c}$ be arbitrary but fixed. Then we have:\\
(i) $I(t \star u) \rightarrow 0$ as $t \rightarrow 0$;\\
(ii) $I(t\star u) \rightarrow-\infty$ as $t \rightarrow+\infty$.

\textbf{Proof.}\ \ (i) By Gagliardo-Nirenberg inequality \eqref{GN}, then we have that
$$
\begin{aligned}
|I(t\star u)|&\leq \frac{a}{2} \int_{\mathbb{R}^{N}}|\nabla(t\star u)|^{2} dx+\frac{b}{4} \left (\int_{\mathbb{R}^{N}}|\nabla(t\star u)|^{2} dx \right )^{2}\\
&+\frac{1}{p}\int_{\mathbb{R}^{N}}|t\star u|^{p} dx+\frac{1}{q} \int_{\mathbb{R}^{N}}h\left|t\star u \right |^{q} dx \\
& \leq \frac{at^{2}}{2}\|\nabla u\|_{2}^{2}+\frac{bt^{4}}{4}\|\nabla u\|_{2}^{4}+\frac{t^{p \gamma_{p}}}{p} \mathcal{C}_{N, p}^{p} c^{\frac{p-p \gamma_{p}}{2}}\|\nabla u\|_{2}^{p \gamma_{p}}+\frac{1}{q} t^{q \gamma_{p}} \mathcal{C}_{N, p}^{q} c^{\frac{q(1-\gamma_{p})}{2}}\|h\|_{\frac{p}{p-q} }\|\nabla u\|_{2}^{q \gamma_{p}} \\
& \rightarrow 0
\end{aligned}
$$
as $t\to 0$, since $p\gamma _{p},q\gamma _{p}>0$.

(ii) Similarly, we have that
$$
\begin{aligned}
I(t\star u) & \leq \frac{at^{2}}{2}\|\nabla u\|_{2}^{2}+\frac{bt^{4}}{4}\|\nabla u\|_{2}^{4}-\frac{1}{p} \int_{\mathbb{R}^{N}}|t \star u|^{p} dx+\frac{1}{q} \int_{\mathbb{R}^{N}}h\left | t\star u \right |^{q} dx \\
& \leq \frac{at^{2}}{2}\|\nabla u\|_{2}^{2}+\frac{bt^{4}}{4}\|\nabla u\|_{2}^{4}-\frac{t^{p\gamma_{p}}}{p} \int_{\mathbb{R}^{N}}|u|^{p} dx+\frac{1}{q} t^{q \gamma_{p}} \mathcal{C}_{N, p}^{q} c^{\frac{q(1-\gamma_{p})}{2}}\|h\|_{\frac{p}{p-q} }\|\nabla u\|_{2}^{q \gamma_{p}} \\
& \rightarrow -\infty
\end{aligned}
$$
as $t \rightarrow+\infty$, since $p \gamma_{p}>4$.\ \ \ $\Box$\\

Again using the Gagliardo-Nirenberg inequality,
\begin{equation}\label{k1}
\begin{aligned}
I(u) &\ge \frac{a}{2}\|\nabla u\|_{2}^{2}+\frac{b}{4}\|\nabla u\|_{2}^{4}-\frac{1}{p} \mathcal{C}_{N, p}^{p}c^{\frac{p-p \gamma_{p}}{2}}\|\nabla u\|_{2}^{p \gamma_{p}}-\frac{1}{q} \mathcal{C}_{N, p}^{q} c^{\frac{q(1-\gamma_{p})}{2}}\| h \| _{\frac{p}{p-q}} \|\nabla u\|_{2}^{q \gamma_{p}}\\
&\ge  \frac{a}{2}\|\nabla u\|_{2}^{2}-\frac{1}{p} \mathcal{C}_{N, p}^{p} c^{\frac{p-p \gamma_{p}}{2}}\|\nabla u\|_{2}^{p \gamma_{p}}-\frac{1}{q} \mathcal{C}_{N, p}^{q} c^{\frac{q(1-\gamma_{p})}{2}}\left \| h \right \| _{\frac{p}{p-q}} \|\nabla u\|_{2}^{q \gamma_{p}}.
\end{aligned}
\end{equation}
To understand the geometry of the functional $I$ on $S_{c}$, it is useful to consider the function $\varphi: \mathbb{R}^{+} \rightarrow \mathbb{R}$ defined by
\begin{equation}\label{k2}
\varphi(t):=\frac{a}{2} t^{2}-\frac{1}{p} \mathcal{C}_{N, p}^{p} c^{\frac{p-p \gamma_{p}}{2}} t^{p \gamma_{p}}-\frac{1}{q} \mathcal{C}_{N, p}^{q} c^{\frac{q(1-\gamma_{p})}{2}}\left \| h \right \| _{\frac{p}{p-q} }t^{q \gamma_{p}}.
\end{equation}
Since $0<q \gamma_{p}<2<p\gamma_{p}$, we have that $\varphi\left(0^{+}\right)=0^{-}$ and $\varphi(+\infty)=-\infty$. The role of assumption \eqref{Lcondi2} is clarified by the following lemma.\\

\textbf{Lemma 3.4.}\ \ Under the assumption $\left(\mathbf{h}_{\mathbf{2}}\right)$, if \eqref{Lcondi1} holds,
then the function $\varphi$ has a local strict minimum at negative level and a global strict
maximum at positive level. Moreover, there exist $0<R_{1}<R_{2},$ both depending on $c$, such that $\varphi(R_{1})=0=\varphi(R_{2})$ and $\varphi(t)>0$ if and only if $t\in (R_{1},R_{2})$.

\textbf{Proof.}\ \ For $t>0$, we see that $\varphi(t)>0$ if and only if
$$
\psi(t)>\frac{1}{q} \mathcal{C}_{N, p}^{q} c^{\frac{q(1-\gamma_{p})}{2}}\left \| h \right \| _{\frac{p}{p-q} },
$$
where
$$
\psi(t):=\frac{a}{2} t^{2-q\gamma_{p}}-\frac{1}{p} \mathcal{C}_{N, p}^{p} c^{\frac{p-p \gamma_{p}}{2}} t^{p \gamma_{p}-q\gamma_{p}}.
$$
Observe that $p \gamma_{p}-q\gamma_{p}>2-q\gamma_{p}>0$, then $\psi$ has a unique critical point $\bar{t}$ on $(0,+\infty)$, which is a global maximum point at positive level. In fact, the expression of $\bar{t}$ is
$$
\bar{t}=\left(\frac{ap\left(2-q\gamma_{p}\right)}{2 \gamma_{p}(p-q) \mathcal{C}_{N, p}^{p} c^{\frac{p-p \gamma_{p}}{2}}}\right)^{\frac{1}{p\gamma_{p}-2}},
$$
and the maximum value of $\psi$ is
\begin{equation}\label{k3}
\psi(\bar{t})=\frac{a(p\gamma_{p}-2)}{2 \gamma_{p}(p-q)}\left(\frac{ap\left(2-q\gamma_{p}\right)}{2 \gamma_{p}(p-q) \mathcal{C}_{N, p}^{p}}\right)^{\frac{2-q\gamma_{p}}{p \gamma_{p}-2}} c^{-\frac{p\left(1-\gamma_{p}\right)\left(2-q\gamma_{p}\right)}{2\left(p \gamma_{p}-2\right)}}.
\end{equation}
Therefore, if \eqref{Lcondi2} holds, then $\psi(\bar{t})>\frac{1}{q} \mathcal{C}_{N, p}^{q} c^{\frac{q(1-\gamma_{p})}{2}}\left \| h \right \| _{\frac{p}{p-q}},$ thus the equation $\varphi=0$ has two roots  $R_{1}, R_{2}$ and $\varphi$ is positive on $\left(R_{1}, R_{2}\right)$. Moreover, $\varphi$ has a global maximum point  $t_{2}$. According to the expression of $\varphi$, we can deduce that $\varphi$ also has a local minimum point $t_{1}$ at negative level in $\left(0, R_{1}\right)$.\ \ \ $\Box$\\

Set
$$
\begin{aligned}
A_{\iota}:=\left\{u \in S_{c}:|\nabla u|_{2}<\iota\right\}, \\
I^{k}:=\{u \in S_{c}:I(u)<k\} .
\end{aligned}
$$
By Lemmas 3.3 and 3.4, there exists a $\iota_{1}>0$ small enough, such that
$$
I(u)<\frac{1}{2} \varphi\left(t_{2}\right), \text { for any } u \in A_{\iota_{1}} .
$$
Moreover, $I^{\varphi\left(t_{1}\right)} \subset\left\{|\nabla u|_{2}>R_{2}\right\}$ since $I(u)\geq \varphi\left(|\nabla u|_{2}\right)$.
Now we can get a mountain pass structure of $I$ on manifold $S_{c}$.
\begin{equation}\label{k4}
\Gamma:=\{\gamma \in C([0,1], S_{c}): \gamma(0) \in A_{\iota}, \gamma(1) \in I^{\varphi\left(t_{1}\right)}\},
\end{equation}
and the mountain pass value is
\begin{equation}\label{k5}
m_{h,c}:=\inf _{\gamma \in \Gamma} \max _{t\in[0,1]} I(\gamma(t)).
\end{equation}\\

\textbf{Remark 3.1.}\ \
$$
I_{\infty}\left(v_{c}\right)=m_{c}=\inf_{\gamma \in \Gamma}\max_{t\in[0,1]} I_{\infty}(\gamma(t))
$$
where $v_{c}$ satisfies
$$
\left\{
\begin{aligned}
&-\left(a+b\int_{\mathbb{R}^{N}}\left | \nabla v_{c} \right |^{2}dx\right)\Delta v_{c}+\lambda v_{c}=|v_{c}|^{p-2} v_{c}  \quad \text { in } \mathbb{R}^{N}, \\
&\int_{\mathbb{R}^{N}}\left|v_{c}\right|^{2}dx=c, \quad u \in H^{1}(\mathbb{R}^{N}),
\end{aligned}
\right.
$$
i.e., the solution $v_{c}$ of the problem \eqref{lim-initial} is a mountain pass critical point of $I_{\infty}$ constrained on $S_{c}$.
(see \cite{C1}), it is immediately seen that
\begin{equation}\label{mhm}
m_{h,c}<m_{c}.
\end{equation}

\textbf{Lemma 3.5.}\ \ Under the assumption $\left(\mathbf{h}_{\mathbf{2}}\right)$, suppose that $h$  satisfies \eqref{Lcondi2}, then there exists a $(PS)$ sequence $\{u_{n}\}$ of $\left.I\right|_{S_{c}}$, which satisfies
\begin{equation}\label{w1}
I\left(u_{n}\right) \rightarrow m_{h,c},
\end{equation}
\begin{equation}\label{w2}
I^{\prime} \mid_{S_{c}}(u_{n}) \rightarrow 0, \end{equation}
\begin{equation}\label{w3}
P\left(u_{n}\right) \rightarrow 0,
\end{equation}
as $n \rightarrow \infty$, where
$$
P(u)=a\|\nabla u\|_{2}^{2}+b\|\nabla u\|_{2}^{4}-\gamma_{p}\int_{\mathbb{R}^{N}}|u|^{p} d x-\gamma_{q}\int_{\mathbb{R}^{N}} h \left | u \right |^{q}  dx+\frac{1}{q} \int_{\mathbb{R}^{N}}\langle\nabla h, x\rangle \left | u \right |^{q} dx,
$$
and
\begin{eqnarray}\label{fushu}
\lim_{n \to \infty } \left \| \left ( u_{n} \right )^{-} \right \| =0.
\end{eqnarray}
We remark that \eqref{w2} means that there exists  $\{\lambda_{n}\}_{n \geq 1}$, such that for any $\varphi \in C_{0}^{\infty}(\mathbb{R}^{N})$, there holds
\begin{equation}\label{w4}
I^{\prime}\left(u_{n}\right)[\varphi]+\lambda_{n} \int_{\mathbb{R}^{N}} u_{n} \varphi dx \rightarrow 0,\quad as \ \ n \rightarrow \infty.
\end{equation}
Moreover, $\{u_{n}\}$ is bounded in $H^{1}(\mathbb{R}^{N})$ and the related Lagrange multipliers  $\{\lambda_{n}\}$ in \eqref{w4} is also bounded, up to a subsequence, $\lambda_{n}\to \bar{\lambda }$, with $\bar{\lambda }> 0$.

\textbf{Proof.}\ \ We divide the proof into three steps.

\textbf{Step1}: Existence of the Palais-Smale sequence. The existence of the $(PS)$ sequence that verifies \eqref{w3} and \eqref{fushu} closely follows the arguments in \cite{r18}, where the authors adapt some ideas in \cite{r4}. We recall the main strategy, referring to \cite{r18} for the details. A key tool is to set
$$
\tilde{I}(t, u):=I(t\star u) \quad \text { for all }(t,u) \in \mathbb{R}\times H^{1}(\mathbb{R}^{N}).
$$
The corresponding minimax structure of $\tilde{I}$ on $\mathbb{R} \times S_{c}$ as follows
\begin{equation}\label{k6}
\tilde{\Gamma}:=\{\gamma=\left(\gamma_{1}, \gamma_{2}\right) \in C([0,1], \mathbb{R} \times S_{c}): \gamma(0) \in\left(0, A_{\iota_{1}}\right), \gamma(1) \in(0,I^{\varphi\left(t_{1}\right)})\},
\end{equation}
and it minimax value is
\begin{equation}\label{k7}
\tilde{m}_{h,c}:= \inf _{\gamma \in \tilde{\Gamma}} \max _{t \in[0,1]} \tilde{I}(\gamma(t)).
\end{equation}
It turns out that $\tilde{m}_{h,c}=m_{h,c}$ and that, if $\left(t_{n}, v_{n}\right)_{n}$ is a $(PS)_{c}$ sequence for $\widetilde{I}$ with $t_{n}\to 0$, then  $u_{n}=t_{n} \star v_{n}$ is a $(PS)_{c}$ sequence for $I$. Now, let us consider a sequence $\xi_{n} \in \Gamma$ such that
$$
\begin{aligned}
m_{h, c} \leq \max _{t \in[0,1]} I\left(\xi_{n}(t)\right)<m_{h, c}+\frac{1}{n}
\end{aligned}
$$
We observe that, since $I(u)=I(|u|)$ for every  $u\in H^{1}\left(\mathbb{R}^{N}\right)$, we can take  $\xi_{n}(t) \geq 0$ in $\mathbb{R}^{N}$, for every  $t \in[0,1]$ and  $n \in \mathbb{N}$.
We are in a position to apply Lemma 3.2 to $\tilde{I}$ with
$$
X:=\mathbb{R}\times S_{c}, \quad K:=\{(0, A_{\iota_{1}}), (0,I^{\varphi\left(t_{1}\right)})\}, \quad \mathcal{E}=\widetilde{\Gamma}, \quad E_{n}:=\left\{\left(0,\xi_{n}(t)\right): t \in[0,1]\right\} .
$$
As a consequence, there exist a sequence $(t_{n}, v_{n}) \in \mathbb{R}\times S_{c}$ and $\tilde{c}>0$ such that
\begin{equation}\label{fushuproof}
\begin{aligned}
m_{h,c}-\frac{1}{n}<\tilde{I}\left(t_{n}, v_{n}\right)<m_{h,c}+\frac{1}{n} \\
\min _{t \in[0,1]}\left\|\left(t_{n}, v_{n}\right)-\left(0,\xi_{n}(t)\right)\right\|_{\mathbb{R}\times H^{1}\left(\mathbb{R}^{N}\right) }<\frac{\tilde{c}}{\sqrt{n}} \\
\left\| \nabla_{\mathbb{R}\times S_{c} } \widetilde{I}\left(t_{n}, v_{n}\right)\right\|<\frac{\tilde{c}}{\sqrt{n}} .
\end{aligned}
\end{equation}
Now, we can define
$$
u_{n}=t_{n} \star v_{n}
$$
We observe that, differentiating $\tilde{I}$ with respect to $t$, we get the "almost" Pohozaev identity \eqref{w3}, differentiating with respect to the second variable on the tangent space to $S_{c}$, and by \eqref{fushuproof} and $\xi_{n}(t)\ge0$ we get \eqref{fushu}.

\textbf{Step2}: Boundedness of the $(PS)$ sequence.

By \eqref{w1}, for the $(PS)$ sequence $\{u_{n}\}\subset S_{c}$, there holds
\begin{equation}\label{w5}
\begin{aligned}
m_{h,c}&=I\left(u_{n}\right)+o(1) \\
& =\frac{a}{2}\|\nabla u_{n}\|_{2}^{2}+\frac{b}{4}\|\nabla u_{n}\|_{2}^{4}-\frac{1}{p} \int_{\mathbb{R}^{N}}\left|u_{n}\right|^{p} dx-\frac{1}{q} \int_{\mathbb{R}^{N}} h\left | u_{n} \right |^{q} dx+o(1).
\end{aligned}
\end{equation}
Combining with \eqref{w3},
\begin{equation}\label{w6}
\begin{aligned}
m_{h,c}= & \frac{a(N(p-2)-4)}{2 N(p-2)} \|\nabla u_{n} \|_{2}^{2}+\frac{b(N(p-2)-8)}{4N(p-2)}\|\nabla u_{n} \|_{2}^{4} -\frac{p-q}{q(p-2)} \int_{\mathbb{R}^{N}} h\left | u_{n}\right |^{q} dx \\
& -\frac{2}{qN(p-2)} \int_{\mathbb{R}^{N}}\langle\nabla h, x\rangle \left | u_{n}  \right |^{q} dx+o(1) \\
\ge& \frac{a(N(p-2)-4)}{2 N(p-2)} \|\nabla u_{n} \|_{2}^{2}-\frac{p-q}{q(p-2)} \int_{\mathbb{R}^{N}} h| u_{n}|^{q} dx\\
&-\frac{2}{qN(p-2)} \int_{\mathbb{R}^{N}}\langle\nabla h, x\rangle \left | u_{n}  \right |^{q} dx+o(1)\\
\ge &\frac{a(N(p-2)-4)}{2N(p-2)}\|\nabla u_{n} \|_{2}^{2}-\frac{p-q}{q(p-2)}\mathcal{C}_{N, p}^{q} c^{\frac{q(1-\gamma_{p})}{2}}\|h\|_{\frac{p}{p-q}}\|\nabla u_{n}\|_{2}^{q \gamma_{p}}\\
&-\frac{2}{qN(p-2)}\| \nabla h\cdot x \|_{\frac{2}{2-q}}c^{\frac{q}{2}} +o(1).
\end{aligned}
\end{equation}
Thus $\{u_{n}\}$ is bounded in $H^{1}(\mathbb{R}^{N})$ since $h \in L^{\frac{p}{p-q}}(\mathbb{R}^{N})$ and $\|\nabla h \cdot x\|_{\frac{2}{2-q}}<\infty$.

\textbf{Step3}: Positivity of the Lagrange multiplier.

By taking $u_{n}$ as a test function for \eqref{w4}, we obtain that
$$
o(1)\|u_{n}\|_{H^{1}}=a\|\nabla u_{n}\|_{2}^{2}+b\|\nabla u_{n}\|_{2}^{4}-\| u_{n}  \|_{p}^{p}-\int_{\mathbb{R}^{N}} h \left | u_{n} \right |^{q}  +\lambda_{n}c.
$$
So
$$
\left|\lambda_{n}\right|=\frac{1}{c}\left |  o (1)\left \| u_{n}  \right \|_{H^{1}}-a\|\nabla u_{n}\|_{2}^{2}-b\|\nabla u_{n}\|_{2}^{4}+\| u_{n}\|_{p}^{p}+\int_{\mathbb{R}^{N}} h \left | u_{n} \right |^{q} \right|<+\infty.
$$
Thus the Lagrange multipliers $\{\lambda_{n}\}$ are also bounded.
Next, we show that $\{\lambda_{n}\}$ has a
positive lower bound. In fact, according to \eqref{w3} and \eqref{w4},
\begin{equation}\label{w7}
\begin{aligned}
\lambda_{n}c & =\lambda_{n} \int_{\mathbb{R}^{N}} |u_{n}|^{2} dx \\
& =-a\|\nabla u_{n}\|_{2}^{2}-b\|\nabla u_{n}\|_{2}^{4}+\| u_{n}\|_{p}^{p}+\int_{\mathbb{R}^{N}}h|u_{n}|^{q} dx+o(1) \\
& =\left(1-\gamma_{p}\right) \| u_{n}\|_{p}^{p}+ (1-\gamma_{q})\int_{\mathbb{R}^{N}} h |u_{n}|^{q} dx+\frac{1}{q} \int_{\mathbb{R}^{N}}\langle\nabla h, x\rangle |u_{n}|^{q} dx+o(1).
\end{aligned}
\end{equation}
We also have that
\begin{equation}\label{w8}
\begin{aligned}
m_{h,c}&=\frac{a}{2}\|\nabla u_{n}\|_{2}^{2}+\frac{b}{4}\|\nabla u_{n}\|_{2}^{4}-\frac{1}{p}\|u_{n}\|_{p}^{p}-\frac{1}{q}\int_{\mathbb{R}^{N}}h| u_{n}|^{q}dx+o(1)\\
& =-\frac{b}{4}\|\nabla u_{n}\|_{2}^{4}+\frac{N(p-2)-4}{4p}\|u_{n}\|_{p}^{p}\\
&+\frac{N(q-2)-4}{4q}\int_{\mathbb{R}^{N}}h|u_{n}|^{q}dx-\frac{1}{2q} \int_{\mathbb{R}^{N}}\langle\nabla h, x\rangle |u_{n}|^{q} dx+o(1).
\end{aligned}
\end{equation}
Then combine with the assumption \eqref{Lcondi2}, we have that
\begin{equation}\label{w10}
\begin{aligned}
\lambda_{n}c & +o(1) \\
=&\frac{4p\left(1-\gamma_{p}\right)}{N(p-2)-4} m_{h,c}+\frac{bp\left(1-\gamma_{p}\right)}{N(p-2)-4}\|\nabla u_{n}\|_{2}^{4}+\frac{2p-4}{q(N(p-2)-4)} \int_{\mathbb{R}^{N}}\langle\nabla h, x\rangle |u_{n}|^{q} dx  \\
&+\left (\frac{2q-N(q-2)}{2q}+\frac{(2p-N(p-2))(4-N(q-2))}{2q(N(p-2)-4)}\right )\int_{\mathbb{R}^{N}}h|u_{n}|^{q}dx +o(1)\\
\ge& \frac{4p\left(1-\gamma_{p}\right)}{N(p-2)-4}m_{h,c}-\frac{2p-4}{q(N(p-2)-4)}\| \nabla h\cdot x \|_{\frac{2}{2-q}}c^{\frac{q}{2}}+o(1)\\
\end{aligned}
\end{equation}
since
$$
\begin{aligned}
\|\nabla h \cdot x\|_{\frac{2}{2-q}}<\frac{q(2p-Np+2N)}{p-2}m_{c}c^{-\frac{q}{2}}.
\end{aligned}
$$

Now we prove the convergence of the $(PS)$ sequence $\{u_{n}\}$ and hence we complete the proof of Theorem 1.2.\\

\textbf{Proof of Theorem 1.2.}\ \ Next we prove the existence of solutions of \eqref{initial} with a positive energy level when $2+\frac{8}{N} <p< 2^{*}$. We consider the bounded $(PS)$ sequence $\{u_{n}\}$ given by Lemma 3.5. Then there exists $u \in H^{1}(\mathbb{R}^{N})$ such that $u_{n} \rightharpoonup u$ due to the boundedness of $\{u_{n}\}$. We claim that $u_{n} \rightarrow u$ strongly in $H^{1}(\mathbb{R}^{N}).$

For any $\psi \in H^{1}(\mathbb{R}^{N})$, $\{u_{n}\}$ satisfies
$$
\begin{aligned}
a \int_{\mathbb{R}^{N}} \nabla u_{n} \nabla \psi dx+b\left(\int_{\mathbb{R}^{N}}\left|\nabla u_{n}\right|^{2} dx\right) \int_{\mathbb{R}^{N}} \nabla u_{n} \nabla \psi dx \\
\quad-\int_{\mathbb{R}^{N}}\left|u_{n}\right|^{p-2} u_{n} \psi dx-\int_{\mathbb{R}^{N}} h(x) \left | u_{n} \right |^{q-2}u_{n}   \psi dx \\
\quad=-\lambda_{n} \int_{\mathbb{R}^{N}} u_{n} \psi dx+o(1)\|\psi\|.
\end{aligned}
$$
Using the boundedness of $\{\lambda_{n}\}$ again, we obtain that
$$
\begin{aligned}
a \int_{\mathbb{R}^{N}} \nabla u_{n} \nabla \psi dx+b\left(\int_{\mathbb{R}^{N}}\left|\nabla u_{n}\right|^{2} dx\right) \int_{\mathbb{R}^{N}} \nabla u_{n} \nabla \psi dx \\
\quad-\int_{\mathbb{R}^{N}}\left|u_{n}\right|^{p-2}u_{n} \psi dx-\int_{\mathbb{R}^{N}} h(x) \left | u_{n} \right |^{q-2}u_{n}\psi dx \\
\quad=-\bar{\lambda }\int_{\mathbb{R}^{N}} u_{n} \psi dx+(\bar{\lambda}-\lambda _{n} )\int_{\mathbb{R}^{N}} u_{n} \psi dx+o(1)\|\psi\|.
\end{aligned}
$$
And hence
$$
\begin{aligned}
&a \int_{\mathbb{R}^{N}} \nabla u_{n} \nabla \psi dx+b\left(\int_{\mathbb{R}^{N}}\left|\nabla u_{n}\right|^{2} dx\right) \int_{\mathbb{R}^{N}} \nabla u_{n} \nabla \psi dx \\
\quad&-\int_{\mathbb{R}^{N}}\left|u_{n}\right|^{p-2} u_{n} \psi dx-\int_{\mathbb{R}^{N}} h(x) \left | u_{n} \right |^{q-2}u_{n}\psi dx \\
\quad&=-\bar{\lambda}\int_{\mathbb{R}^{N}} u_{n} \psi dx,
\end{aligned}
$$
which implies that $\{u_{n}\}$ is a $(PS)$ sequence for  $I_{\lambda}$ at level $m_{h,c}+\frac{\lambda}{2} c$, so that we can apply the Splitting Lemma 3.1, getting
$$u_{n}= u+\Sigma_{j=1}^{k} \omega^{j}\left(\cdot-y_{n}^{j}\right)+o(1).
$$
Assume by contradiction that $k\geq 1$,  or, equivalently, that $\|u\|_{2}^{2} <c$. In addition, if $0<\alpha<\beta$, then $m_{\alpha}>m_{\beta}$ and $J_{\infty,0}\left(\omega^{j}\right) \geq m_{\alpha_{j}}$ (see \cite{C7}). Therefore,
\begin{equation}\label{imq}
\begin{aligned}
m_{h,c}+\frac{\lambda}{2} c=J_{h, 0}(u)+\frac{\lambda}{2} \beta+\Sigma_{j=1}^{k} J_{\infty, 0}\left(\omega^{j}\right)+\frac{\lambda}{2} \Sigma_{j=1}^{k} \alpha_{j},
\end{aligned}
\end{equation}
where $\beta:=\|u\|_{2}^{2}$, $\alpha_{j}:=\left\|\omega^{j}\right\|_{2}^{2}$. By \eqref{2fs} we have
$$
c=\beta+\Sigma_{j=1}^{k}\alpha_{j},
$$
thus \eqref{imq}
\begin{equation}\label{fi}
\begin{aligned}
m_{h,c}=J_{h, 0}(u)+\Sigma_{j=1}^{k} J_{\infty, 0}\left(\omega^{j}\right).
\end{aligned}
\end{equation}
Since $J_{h,0}(u), J_{\infty, 0}\left(\omega^{j}\right) \geq m_{c}$, we have $m_{h,c} \geq m_{c}$, which is a contradiction with \eqref{mhm}. Thus  $k=0$. That is $u_{n} \rightarrow u$ strongly in  $H^{1}\left(\mathbb{R}^{N}\right)$ and $u$ is a solution of Eq.\eqref{initial}.\ \ \ $\Box$

\section{Proof of Theorem 1.3}

\ \ \ \ \ \ In this section, we assume that $2+\frac{8}{N} <p< 2^{*}$, $1\le N \le 3$, $\bar{h}(x)=-h(x)\ge 0$ and $\bar{h}(x) \not \equiv 0$. By using a min-max argument, we can find the existence of normalized solutions of Eq.\eqref{initial}. Firstly, we
show that the energy functional corresponding to Eq.\eqref{initial} has a linking geometry. For $s \in \mathbb{R}$ and $u \in H^{1}\left(\mathbb{R}^{N}\right)$, we introduce the scaling
$$
s\star u(x):= e^{\frac{N}{2}s}u\left(e^{s} x\right),
$$
which preserves the $L^{2}$-norm: $\|s * u\|_{2}=\|u\|_{2}$ for all $s \in \mathbb{R}$. For $\mathbb{R}>0$ and $s_{1} <0<s_{2}$, which will be determined later, we set
$$
Q:=B_{R} \times\left[s_{1}, s_{2}\right] \subset \mathbb{R}^{N} \times \mathbb{R}
$$
where $B_{R}=\left\{x \in \mathbb{R}^{N}:|x| \leq R\right\}$ is the closed ball of radius $R$ around 0 in $\mathbb{R}^{N}.$  For $c>0$, define
$$
\Gamma_{c}:=\left\{\gamma: Q \rightarrow S_{c} \mid \gamma \in C\left(\mathbb{R}^{N}\right), \gamma(y, s)=s\star v_{c}(\cdot-y) \text { for all }(y, s) \in \partial Q\right\},
$$
where $v_{c}$ satisfies
$$
\left\{
\begin{aligned}
&-\left(a+b\int_{\mathbb{R}^{N}}\left | \nabla v_{c} \right |^{2}dx\right)\Delta v_{c}+\lambda v_{c}=|v_{c}|^{p-2} v_{c}   \quad \text { in } \mathbb{R}^{N}, \\
&\int_{\mathbb{R}^{N}}\left|v_{c}\right|^{2}dx=c, \quad u \in H^{1}(\mathbb{R}^{N}).
\end{aligned}
\right.
$$

We want to find a solution to Eq.\eqref{initial} in $S_{c}$ whose energy $I$ is given by
$$
L_{h, c}:=\inf _{\gamma \in \Gamma_{c}} \max _{(y, s) \in Q}I(\gamma(y, s)).
$$
To prove that the energy functional $I$ has a linking geometry, it is necessary to find the suitable $R>0$, $s_{1}<0<s_{2}$ such that
$$
\sup _{\gamma \in \Gamma_{c}} \max _{(y, s) \in \partial Q}I(\gamma(y, s))<L_{h,c}
$$
at least for some suitable choice of $Q$. Now we recall the notion of barycentre of a function $u \in H^{1}\left(\mathbb{R}^{N}\right) \backslash\{0\}$ which has been introduced in \cite{p1} and in \cite{p2}. Setting
$$
\nu(u)(x)=\frac{1}{\left|B_{1}(0)\right|} \int_{B_{1}(x)}|u(y)| dy,
$$
we observe that $\nu(u)$ is bounded and continuous, so the function
$$
\hat{u}(x)=\left[\nu(u)(x)-\frac{1}{2} \max \nu(u)\right]^{+}
$$
is well defined, continuous, and has compact support. Therefore we can define $\beta: H^{1}\left(\mathbb{R}^{N}\right) \backslash\{0\} \rightarrow \mathbb{R}^{N}$ as
$$
\beta(u)=\frac{1}{\|\hat{u}\|_{1}} \int_{\mathbb{R}^{N}} \hat{u}(x) x d x .
$$
The map $\beta$ is well defined, because $\hat{u}$ has compact support, and it is not difficult to verify that it enjoys the following properties:\\
(i) $\beta$ is continuous in  $H^{1}\left(\mathbb{R}^{N}\right) \backslash\{0\}$;\\
(ii) if $u$ is a radial function, then $\beta(u)=0$;\\
(iii) $\beta(t u)=\beta(u)$ for all $t \neq 0$ and for all $u \in H^{1}\left(\mathbb{R}^{N}\right) \backslash\{0\}$;\\
(iv) setting $u_{z}(x)=u(x-z)$ for $z \in \mathbb{R}^{N}$ and $u \in H^{1}\left(\mathbb{R}^{N}\right) \backslash\{0\}$ there holds $\beta\left(u_{z}\right)=\beta(u)+z$.\\
Now we define
$$
\begin{aligned}
&\mathcal{D}:=\left\{D \subset S_{c}: D \text { is compact, connected, } s_{1} \star v_{c}, s_{2} \star v_{c} \in D\right\}, \\
&\mathcal{D}_{0}:=\{D \in \mathcal{D}: \beta(u)=0 \text { for all } u \in D\}, \\
&\mathcal{D}_{r}:=\mathcal{D} \cap H_{\text {rad }}^{1}\left(\mathbb{R}^{N}\right),
\end{aligned}
$$
and
$$
\begin{aligned}
w_{c}^{r} & :=\inf _{D \in \mathcal{D}_{r}} \max _{u \in D} I_{\infty}(u) \\
w_{c}^{0} & :=\inf _{D \in \mathcal{D}_{0}} \max _{u \in D} I_{\infty}(u) \\
w_{c} & :=\inf _{D \in \mathcal{D}} \max _{u \in D} I_{\infty}(u) .
\end{aligned}
$$
It has been proved in \cite{C7} that
$$
m_{c}=\inf _{\sigma \in \Sigma_{c}} \max _{t \in[0,1]} I_{\infty}(\sigma(t))
$$
where
$$
\Sigma_{c}=\left\{\sigma \in \mathcal{C}\left([0,1], S_{c}\right): \sigma(0)=s_{1} \star v_{c}, \sigma(1)=s_{2} \star v_{c}\right\} .
$$\\

\textbf{Lemma 4.1.}\ \ $w_{c}^{r}=w_{c}^{0}=w_{c}=m_{c}$.

\textbf{Proof.}\ \ Clearly $\mathcal{D}_{r} \subset \mathcal{D}_{0} \subset \mathcal{D}
$, so that $w_{c}^{r} \geq w_{c}^{0} \geq w_{c}$. It remains to prove that $w_{c} \geq m_{c}$ and  $m_{c} \geq w_{c}^{r}$.

Arguing by contradiction we assume that  $m_{c}>w_{c}$. Then $\max _{u \in D} I_{\infty}(u)<m_{c}$ for some $D \in \mathcal{D}$, hence $\sup _{u \in U_{\delta}(D)} I_{\infty}(u)<m_{c}$ for some $\delta>0$, here  $U_{\delta}(D)$ is the $\delta$-neighborhood of $D$. Observe that $U_{\delta}(D)$ is open and connected, so it is path-connected. Therefore there exists a path $\sigma \in \Sigma_{c}$ such that $\max _{t \in[0,1]} I_{\infty}(\sigma(t))<m_{c}$, a contradiction.

The inequality $m_{c} \geq w_{c}^{r}$ follows from the fact that the set $D:= \left\{s \star v_{c}: s \in\left[s_{1}, s_{2}\right]\right\} \in \mathcal{D}_{r}$ satisfies
$$
\max _{u \in D} I_{\infty}(u)=\max _{s \in\left[s_{1}, s_{2}\right]} I_{\infty}\left(s \star v_{c}\right)=m_{c}.\ \ \  \Box
$$\\

\textbf{Lemma 4.2.}\ \ $L_{c}:=\inf _{D \in \mathcal{D}_{0}} \max _{u \in D} I(u)>m_{c}$.

\textbf{Proof.}\ \  Using $\bar{h}(x) \geq 0$ and Lemma 4.1, we have
\begin{equation}\label{ffmc}
\max _{u \in D} I(u) \geq \max _{u \in D} I_{\infty}(u) \geq w_{c}^{0}=m_{c}, \quad \text {for all} \ D\in \mathcal{D}_{0}.
\end{equation}
Now we argue by contradiction and assume that there exists a sequence $D_{n} \in \mathcal{D}_{0}$ such that
$$
\max _{u \in D_{n}} I(u) \rightarrow m_{c}.
$$
In view of \eqref{ffmc}, we also have
$$
\max _{u \in D_{n}} I_{\infty}(u) \rightarrow m_{c}.
$$
Adapting an argument from \cite[Lemma 2.4]{r4}, we consider the functional
$$
\tilde{I}_{\infty}: H^{1}\left(\mathbb{R}^{N}\right) \times \mathbb{R} \rightarrow \mathbb{R}, \quad \tilde{I}_{\infty}(u, s):=I_{\infty}(s \star u)
$$
constrained to $M:=S_{c} \times \mathbb{R}$. We apply Lemma 3.2 with
$$
K:=\left\{\left(s_{1} \star v_{c}, 0\right),\left(s_{2} \star v_{c}, 0\right)\right\}
$$
and
$$
\mathcal{C}:=\{C \subset M: \ C \text { compact, connected,}\  K \subset C\}.
$$
Observe that
$$
\tilde{w}_{c}:=\inf _{C \in \mathcal{C}} \max _{(u,s) \in C} \tilde{I}_{\infty}(u,s)=w_{c}=m_{c}
$$
because $\mathcal{D} \times\{0\} \subset \mathcal{C}$, hence $w_{c} \geq \tilde{w}_{c}$, and for any $C \in \mathcal{C}$ we have $D:=\{s \star u:(u, s) \in C\} \in \mathcal{D}$ and
$$
\max _{(u, s) \in C} \tilde{I}_{\infty}(u, s)=\max _{(u, s) \in C} I_{\infty}(s \star u)=\max _{v \in D} I_{\infty}(v),
$$
hence $w_{c} \leq \tilde{w}_{c}$. Hence, Lemma 3.2 yields a sequence $\left(u_{n}, s_{n}\right) \in S_{c} \times \mathbb{R}$ such that\\
(1) $\left|\tilde{I}_{\infty}\left(u_{n}, s_{n}\right)-m_{c}\right| \rightarrow 0$ as $n \rightarrow \infty$;\\
(2) $\left\|\nabla_{S_{c} \times \mathbb{R}} \tilde{I}_{\infty}\left(u_{n}, s_{n}\right)\right\| \rightarrow 0$ as  $n \rightarrow \infty$;\\
(3) $\operatorname{dist}\left(\left(u_{n}, s_{n}\right), D_{n} \times\{0\}\right) \rightarrow 0$ as $n \rightarrow \infty.$\\
Then $v_{n}:=s_{n} \star u_{n} \in S_{c}$ is a $(PS)$ sequence for $I_{\infty}$ on $S_{c}$ at $m_{c}$, and there exist Lagrange multipliers  $\lambda_{n} \in \mathbb{R}$ such that
$$
\begin{aligned}
&I_{\infty}\left(v_{n}\right) \rightarrow m_{c}, \\ &a\left\|\nabla v_{n}\right\|_{2}^{2}+b\left\|\nabla v_{n}\right\|_{2}^{4}-\frac{N(p-2)}{2 p}\left\|v_{n}\right\|_{p}^{p} \rightarrow 0, \\
&\left\|I_{\infty}^{\prime}\left(v_{n}\right)+\lambda_{n} G^{\prime}\left(v_{n}\right)\right\|_{\left(H^{1}\left(\mathbb{R}^{N}\right)\right)^{*}} \rightarrow 0, \quad \text { where } G(u)=\frac{1}{2} \int_{\mathbb{R}^{N}}u^{2}dx,
\end{aligned}
$$
as $n \rightarrow \infty$. So, combining those properties we can infer that
$$
\frac{N(p-2)-4}{2N(p-2)} a\left\|\nabla v_{n}\right\|_{2}^{2}+\frac{N(p-2)-8}{4 N(p-2)} b\left\|\nabla v_{n}\right\|_{2}^{4} \rightarrow m_{c}>0, \text { as } n \rightarrow \infty,
$$
and
$$
\begin{aligned}
-\lambda_{n}c&=a\left\|\nabla v_{n}\right\|_{2}^{2}+b\left\|\nabla v_{n}\right\|_{2}^{4}-\left\|v_{n}\right\|_{p}^{p} \\
& =\frac{N(p-2)-2 p}{2 p}\left\|v_{n}\right\|_{p}^{p}=\frac{N(p-2)-2 p}{N(p-2)}\left(a\left\|\nabla v_{n}\right\|_{2}^{2}+b\left\|\nabla v_{n}\right\|_{2}^{4}\right) .
\end{aligned}
$$
Therefore, $\{v_{n}\}$ is bounded in $H^{1}(\mathbb{R}^{N})$ and $\{\lambda_{n}\}$ is bounded in $\mathbb{R}$. We may assume that $v_{n} \rightharpoonup v$ in  $H^{1}(\mathbb{R}^{N})$, $\left\|\nabla v_{n}\right\|_{2}^{2} \rightarrow A^{2}$ , and  $\lambda_{n} \rightarrow \lambda>0$. In fact, $\{v_{n}\}$ is a $(PS)$ sequence for $I_{\infty, \lambda}$ at level $m_{c}+\frac{\lambda}{2} c.$
As a consequence of Lemma 3.1, $v_{n}$ can be rewritten as
$$
v_{n}=v+\sum_{j=1}^{k} w^{j}\left(\cdot-y_{n}^{j}\right)+o(1)
$$
in $H^{1}(\mathbb{R}^{N})$, where $k \geq 0$ and $w^{j} \neq 0,$ $v$ are solutions to
$$
-\left(a+b A^{2}\right) \Delta w+\lambda w=|w|^{p-2}w
$$
and $|y_{n}^{j}| \rightarrow \infty$. Moreover, we get
\begin{equation}\label{xx2}
c=\|v\|_{2}^{2}+\sum_{j=1}^{k}\left\|w^{j}\right\|_{2}^{2}+o(1),
\end{equation}
\begin{equation}\label{xx3}
A^{2}=\|\nabla v\|_{2}^{2}+\sum_{j=1}^{k}\left\|\nabla w^{j}\right\|_{2}^{2},
\end{equation}
$$
I_{\infty, \lambda}\left(v_{n}\right) \rightarrow J_{\infty,\lambda}(v)+\sum_{j=1}^{k} J_{\infty,\lambda}\left(w^{j}\right),
$$
and hence,
$$
m_{c}+\frac{\lambda}{2}c=J_{\infty,0}(v)+\frac{\lambda}{2}\|v\|_{2}^{2}+\sum_{j=1}^{k} J_{\infty,0}\left(w^{j}\right)+\frac{\lambda}{2} \sum_{j=1}^{k}\left\|w^{j}\right\|_{2}^{2}+o(1).
$$
By \eqref{xx2}, we have
$$
m_{c}=J_{\infty,0}(v)+\sum_{j=1}^{k} J_{\infty,0}\left(w^{j}\right)+o(1).
$$
If $v\neq 0$ and $k \geq 1$, we get $A^{2}>\left \|\nabla v \right \|_{2}^{2}$ from \eqref{xx3}, we have
$$
\begin{aligned}
J_{\infty,0}(v)&=\left(\frac{a}{2}+\frac{b A^{2}}{4}\right) \int_{\mathbb{R}^{N}}|\nabla v|^{2} dx-\frac{1}{p} \int_{\mathbb{R}^{N}}|v|^{p} dx.\\
&> \frac{a}{2}\int_{\mathbb{R}^{N}}|\nabla v|^{2} dx+\frac{b}{4}\left ( \int_{\mathbb{R}^{N}}|\nabla v|^{2} dx \right )^{2} dx-\frac{1}{p} \int_{\mathbb{R}^{N}}|v|^{p} dx\\
&=I_{\infty}(v)\\
&\ge m_{\left \| v \right \|_{2}^{2}}\ge m_{c}.
\end{aligned}
$$
Similarly, we have $J_{\infty,0}\left(w^{j}\right)\ge m_{c}$. Thus,
$$
m_{c}+o(1)\geq(k+1) m_{c}+o(1),
$$
we get a contradiction. Therefore, $k=1$ and $v=0$, or $k=0$ and $v \neq 0$.
If $k=1$ and $v=0$, then $v_{n}\left(\cdot+y_{n}^{1}\right)+o(1)=w^{1}$. On the other hand, due to point (3) that $\operatorname{dist}\left(\left(u_{n}, s_{n}\right), D_{n} \times\{0\}\right) \rightarrow 0$, we obtain
$$
\beta\left(w^{1}\right)=\beta\left(v_{n}\left(\cdot+y_{n}^{1}\right)\right)+o(1)=y_{n}^{1}+o(1),
$$
which contradicts the fact $\beta$ is continuous and $\left|y_{n}^{1}\right| \rightarrow \infty$.

If $k=0$ and $v \neq 0$, then $v_{n} \rightarrow v$ in $H^{1}\left(\mathbb{R}^{N}\right)$. Using again point (3), we also have $\beta(v)=0$.
Hence, by the uniqueness, $v_{n} \rightarrow \pm v_{c}$ in $H^{1}\left(\mathbb{R}^{N}\right)$. This implies
$$
I\left(v_{n}\right)=I_{\infty}\left(v_{n}\right)+\frac{1}{q} \int_{\mathbb{R}^{N}} \bar{h}(x) |v_{n}|^{q} dx \rightarrow m_{c}+\frac{1}{q} \int_{\mathbb{R}^{N}} \bar{h}(x) |v_{c}|^{q} dx>m_{c},
$$
which is a contradiction.\ \ \ $\Box$ \\

\textbf{Lemma 4.3.}\ \ For any $c >0$, then $L_{h,c} \ge L_{c}$ holds.

\textbf{Proof.}\ \  Similar to \cite[Proposition 3.5]{r18},  we omit it.\ \ \ $\Box$ \\

\textbf{Lemma 4.4.}\ \ For any $c>0$ and for any $\varepsilon >0$, there exist $\bar{R}>0$ and  $\bar{s}_{1}<0<\bar{s}_{2}$ such that for $Q=B_{R} \times\left[s_{1}, s_{2}\right]$ with $R \geq \bar{R}, s_{1} \leq \bar{s}_{1}, s_{2} \geq \bar{s}_{2}$ the following holds:
$$
\max _{(y, s) \in \partial Q} I\left(s \star v_{c}(\cdot-y)\right)<m_{c}+ \varepsilon.
$$

\textbf{Proof.}\ \ We have
$$
I\left(s \star v_{c}(\cdot-y)\right)=I_{\infty}\left(s \star v_{c}\right)+\frac{e^{\frac{qsN}{2}}}{q} \int_{\mathbb{R}^{N}}\bar{h}(x)v_{c}\left(e^{s}(x-y)\right)^{q} dx
$$
and
$$
\begin{aligned}
I_{\infty}\left(s \star v_{c}\right) & =\frac{e^{2 s}}{2} \int_{\mathbb{R}^{N}}\left|\nabla v_{c}\right|^{2}dx+\frac{e^{4s}}{4} \left(\int_{\mathbb{R}^{N}}\left|\nabla v_{c}\right|^{2}dx \right)^{2}  -\frac{e^{\frac{N}{2}(p-2) s}}{p} \int_{\mathbb{R}^{N}} |v_{c}|^{p} dx \\
& =\left\{\begin{array}{ll}
O\left(-e^{\frac{N}{2}(p-2)s}\right) \rightarrow-\infty & \text { as } s \rightarrow \infty, \\
O\left(e^{2s}\right) \rightarrow 0 & \text { as } s \rightarrow-\infty .
\end{array}\right.
\end{aligned}
$$
Moreover, there holds
$$
\begin{aligned}
\frac{e^{\frac{qsN}{2}}}{q}\int_{\mathbb{R}^{N}} \bar{h}(x) v_{c}^{q} \left(e^{s}(x-y)\right)dx & \le \frac{e^{\frac{qsN}{2}}}{q} \left (\int_{\mathbb{R}^{N}} \bar{h}^{\frac{2}{2-q}}dx\right)^{\frac{2-q}{2}} \left ( \int_{\mathbb{R}^{N}}v_{c}^{2}\left(e^{s}(x-y)\right)dx \right )^{\frac{q}{2} } \\
&=\frac{1}{q}\| \bar{h}\|_{\frac{2}{2-q}}c^{\frac{q}{2}}\\
\end{aligned}
$$
because $\bar{h}(x)$ satisfying \eqref{th13}, thus for all $s\in \mathbb{R}$, we have
$$
\frac{e^{\frac{qsN}{2} }}{q}\int_{\mathbb{R}^{N}} \bar{h}(x) v_{c}^{q} \left(e^{s}(x-y)\right)dx<m_{c}.
$$
As a consequence we deduce
$$
\max _{y \in B_{R}, s\in \{s_{1}, s_{2}\}} I\left(s \star v_{c}(\cdot-y)\right)<m_{c}++o(1)
$$
provided $s_{1}<0$ is small enough and $s_{2}>0$ is large enough. Moreover, for $\left |y\right |=R$ large enough and $s \in\left[s_{1}, s_{2}\right]$, we choose $\alpha \in(0,1)$ such that $\alpha\left(1+e^{-s_{1}}\right)<1$ so that we have
$$
\begin{aligned}
&\frac{e^{\frac{qsN}{2}}}{q}\int_{\mathbb{R}^{N}}\bar{h}(x) v_{c}^{q} \left(e^{s}(x-y)\right)dx \\
&\leq \frac{e^{\frac{qsN}{2}}}{q} \int_{|x|>\alpha R} \bar{h}(x) v_{c}^{q} \left(e^{s}(x-y)\right)dx+\frac{e^{\frac{qsN}{2} }}{q}\int_{|x-y|>\alpha Re^{-s}} \bar{h}(x) v_{c}^{q} \left(e^{s}(x-y)\right)dx.
\end{aligned}
$$
The first integral is bounded by
$$
\begin{aligned}
\frac{e^{\frac{qsN}{2}}}{q} \int_{|x|>\alpha R} \bar{h}(x) v_{c}^{q} \left(e^{s}(x-y)\right)dx &\le \frac{e^{\frac{qsN}{2}}}{q} \left ( \int_{|x|>\alpha R} \bar{h}^{\frac{2}{2-q}}dx \right )^{\frac{2-q}{2}}\left (\int_{|x|>\alpha R} v_{c}^{2} \left(e^{s}(x-y)\right)dx\right)^{\frac{q}{2}}\\
&\le \frac{1}{q} \left ( \int_{|x|>\alpha R} \bar{h}^{\frac{2}{2-q} } dx \right )^{\frac{2-q}{2}}\left (\int_{\mathbb{R}^{N}} v_{c}^{2}dx\right)^{\frac{q}{2}}\to 0
\end{aligned}
$$
as $R \rightarrow \infty$ and
$$
\begin{aligned}
\frac{e^{\frac{qsN}{2}}}{q}\int_{|x-y|>\alpha R e^{-s}} \bar{h}(x) v_{c}^{q} \left(e^{s}(x-y)\right)dx&\le \frac{1}{q} \left ( \int_{|x-y|>\alpha R e^{-s}} \bar{h}^{\frac{2}{2-q} } dx \right )^{\frac{2-q}{2}}\left (\int_{|\xi|>\alpha R} v_{c}^{2} \left(\xi\right)d\xi\right)^{\frac{q}{2}}\\
&\le \frac{1}{q} \left ( \int_{\mathbb{R}^{N}} \bar{h}^{\frac{2}{2-q}} dx \right )^{\frac{2-q}{2}}\left (\int_{|\xi|>\alpha R} v_{c}^{2}dx\right)^{\frac{q}{2}}\to 0
\end{aligned}
$$
as $R \rightarrow \infty$, which concludes the proof.\ \ \ $\Box$\\

By Lemma 4.3 and 4.4, we may choose $R>0$ and $s_{1}<0<s_{2}$ such that
$$
\max _{(y, s) \in \partial Q} I\left(s\star v_{c}(\cdot-y)\right)<L_{h,c}.
$$
Therefore, $I$ has a linking geometry and there exists a $(PS)$ sequence at the level $L_{h,c}$. In order to estimate $L_{h,c}$, we have the following Lemma.\\

\textbf{Lemma 4.5.}\ \ If $\left|s_{1}\right|, s_{2}$ are large enough, then
$$
L_{h,c} < 2m_{c}.
$$

\textbf{Proof.}\ \ This follows from
$$
\begin{aligned}
L_{h,c} & \leq \max _{(y, s) \in Q}\left\{I_{\infty}\left(s\star v_{c}(\cdot-y)\right)+\frac{1}{q}\int_{\mathbb{R}^{N}} \bar{h}(x)\left(s\star v_{c}\right)^{q}(x-y) dx\right\} \\
& \leq m_{c}+\frac{1}{q}\left | \bar{h} \right |_{\frac{2}{2-q}}c^{\frac{q}{2}}\\
&< 2m_{c}
\end{aligned}
$$
provided $\left|s_{1}\right|, s_{2}$ are large enough.\ \ \ $\Box$\\

By the Lemma 4.3 and Lemma 4.5, we can get
$$
m_{c}<L_{h,c}<2m_{c}.
$$

Next, we construct a bounded $(PS)$ sequence
of I at $L_{h,c}$ by adopting the approach from \cite{r4} and Lemma 3.2. We define a auxiliary $\mathcal{C}^{1}$ functional
$$
\tilde{I}(u, s):=I(s \star u) \text { for all }(u, s) \in H^{1}\left(\mathbb{R}^{N}\right) \times \mathbb{R},
$$
$$
\tilde{\Gamma}_{c}:=\left\{\tilde{\gamma}: Q \rightarrow S_{c} \mid \tilde{\gamma} \in C\left(\mathbb{R}^{N}\right), \tilde{\gamma}(y, s)=s\star v_{c}(\cdot-y) \text { for all }(y,s) \in \partial Q\right\},
$$
and
$$
\tilde{L}_{h, c}:=\inf _{\tilde{\gamma} \in \tilde{\Gamma_{c}}} \max _{(y, s) \in Q}\tilde{I}(\tilde{\gamma}(y,s)).
$$\\

\textbf{Lemma 4.6.}\ \ (1) $\tilde{L}_{h, c}=L_{h,c}.$ \\
(2) If $\left(u_{n}, s_{n}\right)$ is a $(PS)$ sequence for $\tilde{I}$ at level $\tilde{L}_{h, c}$ and $s_{n} \rightarrow 0$ , then $\left(s_{n} \star u_{n}\right)_{n}$ is a $(PS)$ sequence for $I$ at level $L_{h, c}$.

\textbf{Proof.}\ \ The proof is similar to that of \cite{r4} and is omitted.\ \ \ $\Box$\\

\textbf{Lemma 4.7.}\ \ Let $\tilde{g}_{n} \in \tilde{\Gamma}_{c}$ be a sequence such that
$$
\max _{(y, s) \in Q} \tilde{I}\left(\tilde{g}_{n}(y, s)\right) \leq L_{h,c}+\frac{1}{n} .
$$
Then, there exists a sequence $\left(u_{n}, s_{n}\right) \in S_{c} \times \mathbb{R}$ and $\tilde{c}>0$ such that
$$
\begin{aligned}
L_{h,c}-\frac{1}{n} \leq \tilde{I}\left(u_{n}, s_{n}\right) \leq L_{h,c}+\frac{1}{n} \\
\min _{(y, s) \in Q}\left\|\left(u_{n}, s_{n}\right)-\tilde{g}_{n}(y, s)\right\|_{H^{1}\left(\mathbb{R}^{N}\right) \times \mathbb{R}} \leq \frac{\tilde{c}}{\sqrt{n}}\\
\left\|\nabla_{S_{c} \times \mathbb{R}} \tilde{I}\left(u_{n}, s_{n}\right)\right\| \leq \frac{\tilde{c}}{\sqrt{n}} .
\end{aligned}
$$
The last inequality means:
$$
\left|D \tilde{I}\left(u_{n},s_{n}\right)[(z, s)]\right| \leq \frac{\tilde{c}}{\sqrt{n}}\left(\|z\|_{H^{1}\left(\mathbb{R}^{N}\right)}+|s|\right)
$$
for all
$$
(z, s) \in\left\{(z,s) \in H^{1}\left(\mathbb{R}^{N}\right) \times \mathbb{R}: \int_{\mathbb{R}^{N}} z u_{n} dx=0\right\}
$$

\textbf{Proof.}\ \ Apply Lemma 3.2 to $\tilde{I}$ with
$$
X:=S_{c} \times \mathbb{R}, \quad K:=\left\{\left ( s\star v_{c}(\cdot -y),0  \right):(y,s) \in \partial Q \right\}, \quad \mathcal{E}=\widetilde{\Gamma}_{c}, \quad E_{n}:=\left\{\tilde{g}_{n}(y,s): (y,s)\in Q\right\}.
$$\\

\textbf{Lemma 4.8.}\ \ Under the assumption $\left(\mathbf{h}_{\mathbf{3}}\right)$, then there exists a bounded $(PS)$ sequence $\{v_{n}\}$ of $\left.I\right|_{S_c}$, which satisfies
\begin{equation}\label{z1}
I\left(v_{n}\right) \rightarrow L_{h,c}.
\end{equation}
\begin{equation}\label{z2}
I^{\prime} \mid_{S_{c}}(v_{n}) \rightarrow 0, \end{equation}
\begin{equation}\label{z3}
P\left(v_{n}\right) \rightarrow 0,
\end{equation}
as $n \rightarrow \infty$, where
$$
P(u)=a\|\nabla u\|_{2}^{2}+b\|\nabla u\|_{2}^{4}-\gamma_{p}\int_{\mathbb{R}^{N}}|u|^{p} dx+\gamma_{q}\int_{\mathbb{R}^{N}}\bar{h}\left|u \right |^{q} dx-\frac{1}{q} \int_{\mathbb{R}^{N}}\langle\nabla \bar{h}, x\rangle \left| u\right |^{q} dx,
$$
\begin{eqnarray}\label{zfushu}
\lim_{n \to \infty } \left \| \left (v_{n} \right )^{-}\right \| =0.
\end{eqnarray}
Moreover, the sequence of Lagrange multipliers satisfies, up to subsequence $\lambda_{n}\to \lambda>0$.

\textbf{Proof.}\ \ First, The existence of the $(PS)$ sequence that verifies \eqref{z3} and \eqref{zfushu} closely follows the arguments in Lemma 3.5. The proof is omitted.

Next, we prove $\{v_{n}\}$ is bounded in $H^{1}(\mathbb{R}^{N}).$ By \eqref{z1}, for the $(PS)$ sequence $\left\{v_{n}\right\}\subset S_{c}$, there holds
\begin{equation}\label{z5}
\begin{aligned}
L_{h,c}&=I\left(v_{n}\right)+o(1) \\
&=\frac{a}{2}\|\nabla v_{n}\|_{2}^{2}+\frac{b}{4}\|\nabla v_{n}\|_{2}^{4}-\frac{1}{p}\int_{\mathbb{R}^{N}}\left|v_{n}\right|^{p} dx+\frac{1}{q}\int_{\mathbb{R}^{N}} \bar{h}\left |v_{n} \right |^{q}dx+o(1).
\end{aligned}
\end{equation}
Combining with \eqref{z3},
\begin{equation}\label{z6}
\begin{aligned}
L_{h,c}= & \frac{a(N(p-2)-4)}{2 N(p-2)}\|\nabla v_{n} \|_{2}^{2}+\frac{b(N(p-2)-8)}{4N(p-2)}\|\nabla v_{n}\|_{2}^{4}+\frac{p-q}{q(p-2)} \int_{\mathbb{R}^{N}} \bar{h}\left |v_{n}\right |^{q} dx \\
&+\frac{2}{qN(p-2)}\int_{\mathbb{R}^{N}}\langle\nabla \bar{h}, x\rangle \left |v_{n}\right |^{q} dx+o(1) \\
\ge& \frac{a(N(p-2)-4)}{2N(p-2)}\|\nabla v_{n} \|_{2}^{2}+\frac{2}{qN(p-2)} \int_{\mathbb{R}^{N}}\langle\nabla \bar{h}, x\rangle \left |v_{n}\right |^{q} dx+o(1)\\
\ge &\frac{a(N(p-2)-4)}{2N(p-2)}\|\nabla v_{n}\|_{2}^{2}-\frac{2}{qN(p-2)}\| \nabla \bar{h}\cdot x \|_{\frac{2}{2-q}}c^{\frac{q}{2}} +o(1).
\end{aligned}
\end{equation}
Thus $\{v_{n}\}$ is bounded in $H^{1}(\mathbb{R}^{N})$ since $\|\nabla \bar{h} \cdot x\|_{\frac{2}{2-q}}<\infty$.

Then, we prove the positivity of the Lagrange multiplier. In the same way as lemma 3.5. By \eqref{z2}, we obtain that
$$
\left|\lambda_{n}\right|=\frac{1}{c}\left |o (1)\left \| v_{n}\right \|_{H^{1}}-a\|\nabla v_{n}\|_{2}^{2}-b\|\nabla v_{n}\|_{2}^{4}+\|v_{n}\|_{p}^{p}-\int_{\mathbb{R}^{N}} \bar{h}\left|v_{n} \right |^{q} \right|<+\infty.
$$
Thus the Lagrange multipliers $\{\lambda_{n}\}$ are also bounded.
In fact, according to \eqref{z2} and \eqref{z3}, we have that
\begin{equation}\label{z10}
\begin{aligned}
&\lambda_{n}c+o(1) \\
&=\frac{4p\left(1-\gamma_{p}\right)}{N(p-2)-4}L_{h,c}+\frac{bp\left(1-\gamma_{p}\right)}{N(p-2)-4}\|\nabla v_{n}\|_{2}^{4}-\frac{4(p-q)}{q\left(N(p-2)-4 \right)}\int_{\mathbb{R}^{N}}\bar{h} v_{n}^{q} dx \\
&-\frac{2p-4}{q(N(p-2)-4)} \int_{\mathbb{R}^{N}}\langle\nabla \bar{h}, x\rangle v_{n}^{q} dx +o(1)\\
&\ge \frac{4p\left(1-\gamma_{p}\right)}{N(p-2)-4}m_{c}-\frac{4(p-q)}{q\left ( N(p-2)-4 \right)}\|\bar{h}\|_{\frac{2}{2-q}}c^{\frac{q}{2}}-\frac{2p-4}{q(N(p-2)-4)} \Upsilon \| \bar{h}\|_{\frac{2}{2-q}}c^{\frac{q}{2}}\\
&= \frac{2}{N(p-2)-4}\left ( 2p\left(1-\gamma_{p}\right)m_{c}-\frac{2(p-q)}{q}\| \bar{h}\|_{\frac{2}{2-q}}c^{\frac{q}{2}}-\frac{p-2}{q}\Upsilon \|\bar{h}\|_{\frac{2}{2-q}}c^{\frac{q}{2}} \right )
\end{aligned}
\end{equation}
thus $\lambda>0$ provided
$$
\frac{2(p-q)}{q}\| \bar{h}\|_{\frac{2}{2-q}}c^{\frac{q}{2}}+\frac{p-2}{q}\Upsilon \|\bar{h}\|_{\frac{2}{2-q}}c^{\frac{q}{2}}<2p\left(1-\gamma_{p}\right)m_{c}.
$$
so
$$
\| \bar{h}\|_{\frac{2}{2-q}}<\frac{2p\left(1-\gamma_{p}\right)}{2(p-q)+(p-2)\Upsilon}\cdot \frac{qm_{c}}{c^{\frac{q}{2}}},
$$
which is given in \eqref{th13}.\ \ \ $\Box$\\

\textbf{Proof of Theorem 1.3.}\ \ Since $\{v_{n}\}$ is bounded, after passing to a subsequence it converges weakly in $H^{1}(\mathbb{R}^{N})$ to $v\in H^{1}(\mathbb{R}^{N})$. By \eqref{zfushu} and weak convergence, $v$ is a nonnegative weak solution of
\begin{equation}\label{vsolt}
-\left(a+bA^{2}\right)\triangle v+\lambda v+\bar{h}(x)|v|^{q-2} v=|v|^{p-2} v
\end{equation}
such that $\beta:=\|v\|_{2}^{2} \leq c$,  where $A^{2}:=\lim _{n \rightarrow \infty}\left\|\nabla v_{n}\right\|_{2}^{2}$. We note that $\{v_{n}\}$ is a bounded $(PS)$ sequence of $I_{\lambda}$ at level $L_{h, c}+\frac{\lambda}{2}c$, therefore, by Lemma 3.1, there exists an integer $k\geq 0$, $k$ non-trivial solutions $w^{1}, w^{2}, \ldots, w^{k}$ to the equation
$$
-\left(a+b A^{2}\right) \triangle w+\lambda w=|w|^{p-2} w
$$
and $k$ sequences $\{y_{n}^{j}\} \in H^{1}(\mathbb{R}^{N}), 1 \leq j \leq k$, such that $|y_{n}^{j}| \rightarrow \infty$ as $n \rightarrow \infty$.

Moreover, we have
\begin{equation}\label{flyy1}
\begin{aligned}
&v_{n}-\sum_{j=1}^{k} w^{j}\left(\cdot-y_{n}^{j}\right) \rightarrow v \ \text { in } \ H^{1}\left(\mathbb{R}^{N}\right), \\
&\left\|v_{n}\right\|_{2}^{2} \rightarrow\|v\|_{2}^{2}+\sum_{j=1}^{k}\left\|w^{j}\right\|_{2}^{2}, \quad A^{2}=\|\nabla v\|_{2}^{2}+\sum_{j=1}^{k}\left\|\nabla w^{j}\right\|_{2}^{2},
\end{aligned}
\end{equation}
and
\begin{equation}\label{flyy2}
I_{\lambda}\left(v_{n}\right) \rightarrow J_{h,\lambda}(v)+\sum_{j=1}^{k} J_{\infty,\lambda}\left(w^{j}\right)
\end{equation}
as $n \rightarrow \infty$. It remains to show $k=0$, so that $v_{n} \rightarrow v$ strongly in $ H^{1}(\mathbb{R}^{N})$ and we are done. Thus, by contradiction, we can assume that $k \geq 1$, or equivalently $\beta<c$.

First we exclude the case $v=0$. In fact, if $v=0$ and $ k=1$, we have $w^{1}>0$ and $\left\|w^{1}\right\|_{2}^{2}=c$ and $\left\|\nabla w^{1}\right\|_{2}^{2}=A^{2}$ so that \eqref{flyy2} would give $L_{h, c}=m_{c}$, which is not possible due to Lemma 5.3. On the other hand, if $k \geq 2$, we get  $J_{\infty,0}\left(w^{j}\right) \geq m_{\alpha_{j}}\left(\alpha_{j}:=\left\|w^{j}\right\|_{2}^{2}\right)$ and $m_{\alpha_{j}}>m_{c}$, thus $L_{h,c}> 2m_{c}$, which contradicts with Lemma 4.5.

Therefore from now on we will assume $v\neq0$ and $k\ge 1$. From \eqref{flyy2} and $I(v_{n})\to L_{h,c}$, we deduce
$$
L_{h,c}+\frac{\lambda}{2} c=J_{h,0}(v)+\frac{\lambda}{2} \beta+\sum_{j=1}^{k} J_{\infty,0}\left(w^{j}\right)+\sum_{j=1}^{k}\frac{\lambda}{2} \alpha_{j}
$$
Using \eqref{flyy1}, we have
$$
L_{h,c}= J_{h,0}(v)+\sum_{j=1}^{k} J_{\infty,0}\left(w^{j}\right) .
$$
Then, from $A^{2}>\left \|\nabla v \right \|_{2}^{2}$ and $\bar{h}(x) \ge 0$, we have
$$
\begin{aligned}
J_{h,0}(v)&=\left(\frac{a}{2}+\frac{b A^{2}}{4}\right) \int_{\mathbb{R}^{N}}|\nabla v|^{2} dx+\frac{1}{q}\int _{\mathbb{R}^{N}}\bar{h}\left |v \right |^{q} dx -\frac{1}{p} \int_{\mathbb{R}^{N}}|v|^{p} dx.\\
&\ge \frac{a}{2}\int_{\mathbb{R}^{N}}|\nabla v|^{2} dx+\frac{b}{4}\left ( \int_{\mathbb{R}^{N}}|\nabla v|^{2} dx \right )^{2} dx-\frac{1}{p} \int_{\mathbb{R}^{N}}|v|^{p} dx\\
&=I_{\infty}(v)\\
&\ge m_{\left \| v \right \|_{2}^{2}}\ge m_{c}.
\end{aligned}
$$
Similarly, we have $J_{\infty,0}\left(w^{j}\right)\ge m_{c}$. Thus,
$$
m_{c}+o(1)\geq(k+1) m_{c}+o(1),
$$
we get a contradiction. Thus $k=0$ and $\{v_{n}\}$ converges strongly to $v$ in $H^{1}(\mathbb{R}^{N})$.\ \ \ $\Box$\\

\end{document}